\documentclass[10pt]{amsart}
\usepackage{amscd}
\usepackage{amsmath}
\usepackage{amsfonts}
\usepackage{amssymb}
\usepackage{amsthm}
\usepackage{enumerate}
\usepackage[T1]{fontenc}
\usepackage{hyperref}
\usepackage{mathrsfs}
\usepackage{stackrel}
\usepackage[all]{xy}
\usepackage{yhmath}

\newtheorem{thm}{Theorem}[section]
\newtheorem{lem}[thm]{Lemma}
\newtheorem{defn}[thm]{Definition}
\newtheorem{prop}[thm]{Proposition}
\newtheorem{cor}[thm]{Corollary}

\newtheorem*{rmk}{Remark}
\newtheorem*{rmks}{Remarks}

\DeclareMathOperator{\Sp}{Sp}

\newcommand{\A}{\mathcal{A}}
\newcommand{\B}{\mathcal{B}}

\newcommand{\D}{\mathcal{D}}

\newcommand{\F}{\mathcal{F}}
\newcommand{\G}{\mathcal{G}}
\newcommand{\I}{\mathcal{I}}

\renewcommand{\L}{\mathcal{L}}

\newcommand{\M}{\mathcal{M}}

\newcommand{\N}{\mathcal{N}}

\newcommand{\R}{\mathcal{R}}

\newcommand{\h}[1]{\widehat{#1}}

\newcommand{\w}[1]{\wideparen{#1}}

\begin{document}
	\title[Auslander regularity and almost mathematics]{Auslander regularity of $p$-adic Banach algebras via almost mathematics}
	\author{Andreas Bode}
	\maketitle
	\begin{abstract}
		We discuss an `almost' version of Auslander regularity and use it to prove the Auslander regularity of various Banach algebras over non-discretely valued fields appearing naturally in $p$-adic locally analytic representation theory: completed Weyl algebras, the completed enveloping algebra of a Lie algebra, and the Banach completion of the distribution algebra $D(G, K)$ for a compact $p$-adic Lie group $G$. 
	\end{abstract}
	\tableofcontents
	
\section{Introduction}
Let $K$ be a complete nonarchimedean field (i.e. a field which is complete with respect to some nonarchimedean norm) of mixed characteristic $(0, p)$. Let $R$ denote the valuation ring of $K$, consisting of all elements of norm less than or equal to $1$, with maximal ideal $\mathfrak{m}$ and $k=R/\mathfrak{m}$. Let $\pi\in \mathfrak{m}$ with $\pi\neq 0$. 

This paper is concerned with the algebraic properties of various $K$-Banach algebras arising from $p$-adic representation theory. To be precise, we are interested in the following constructions:
\begin{enumerate}[(i)]
	\item Let $\mathfrak{g}$ be a finite-dimensional Lie algebra over $K$. Fix a $K$-basis $x_1, \hdots, x_m$ and assume after appropriate rescaling that the $R$-span $\mathfrak{g}_R=\sum R\cdot x_i$ is closed under the Lie bracket. Then we set
	\begin{equation*}
		U_n:=\h{U_R(\pi^n\mathfrak{g}_R)}\otimes_R K,
	\end{equation*}
	which is isomorphic to the Banach completion of the universal enveloping algebra $U(\mathfrak{g})$ with respect to the semi-norm induced by the lattice $U_R(\pi^n\mathfrak{g}_R)$. If we take the limit over $n$, we recover the Arens-Michael envelope (hyperenveloping algebra) $\w{U(\mathfrak{g})}$, studied e.g. in \cite{Schmidt}.
	\item Let $m\geq 1$ and let $T_m=K\langle x_1, \hdots, x_m\rangle$ be the Tate algebra in $m$ variables. Then we let
	\begin{equation*}
		\h{A_{m, n}}:=K\langle x_1, \hdots, x_m, \pi^n\frac{\mathrm{d}}{\mathrm{d}x_1}, \hdots, \pi^n\frac{\mathrm{d}}{\mathrm{d}x_m}\rangle
	\end{equation*}
	denote the completed Weyl algebra (with deformation $\pi^n$) -- these are the Tate-Weyl algebras appearing in \cite{AWannals}. We also remark that if we set $X=\Sp T_m$ to denote the $m$-dimensional rigid analytic polydisk, then the Fr\'echet--Stein algebra of analytic differential operators (studied e.g. in \cite{DcapOne}, \cite{SixOp}) on $X$ is precisely given by
	\begin{equation*}
		\w{\D}_X(X)=\varprojlim_n \h{A_{m, n}}.
	\end{equation*}
	\item Let $L$ be a finite extension of $\mathbb{Q}_p$ and let $G$ be a compact, locally $L$-analytic group. Schneider--Teitelbaum defined in \cite{ST} the distribution algebra $D(G, K)$ of locally $L$-analytic distributions with values in $K$, whose module theory encodes the locally analytic representation theory of $G$. The algebra $D(G, K)$ can be written as an inverse limit of Noetherian Banach algebras $D_r(G, K)$, with $1/p\leq r<1$, tending to $1$ from below. 
\end{enumerate}

We prove the following theorem.
\begin{thm}
	\label{MainThm1}
	\leavevmode
	\begin{enumerate}[(i)]
		\item The algebras $U_n$ are Auslander regular for $n\geq 1$, with $\mathrm{gl.dim.}U_n=\mathrm{dim}_K\mathfrak{g}$.
		\item The algebras $\h{A_{m, n}}$ are Auslander regular for $n\geq 1$, with $\mathrm{gl.dim.}\h{A_{m,n}}=m$.
		\item For each $1/p\leq s<1$, there exists $r<1$ with $s\leq r$ such that $D_r(G, K)$ is Auslander regular, with $\mathrm{gl.dim.}D_r(G, K)=\mathrm{dim}_L G$.
	\end{enumerate}
\end{thm}

The definition of Auslander regularity will be recalled in subsection 3.1 -- it is a regularity condition for possibly non-commutative rings which among other things ensures a well-behaved duality and dimension theory (compare \cite[chapter III, section 4]{LVO}). A commutative Noetherian ring is Auslander regular if and only if it has finite global dimension (the non-trivial implication being \cite[chapter III, Example 2.4.3]{LVO}) if and only if it is regular of finite Krull dimension (see \cite[tag 00OE]{stacksproj}), and  in this case the dimension theory alluded to above coincides with the usual notion of dimension of the support of a given module (cf. \cite[chapter III, Theorem 4.1.7]{LVO}). 

In this spirit, Theorem \ref{MainThm1} allows us to apply the framework from \cite[section 8]{ST} to the Fr\'echet--Stein algebras $\w{U(\mathfrak{g})}$, $\w{\D_X}(X)$ and $D(G, K)$ mentioned above, establishing a well-behaved dimension theory in these cases.

Note that in the case of a discretely valued base field, the Auslander regularity of $U_n$ and $\h{A_{m, n}}$ is a result of Ardakov--Wadsley \cite[Theorem 3.3]{AWannals}, while Schmidt proved the analogous result for norm-completed distribution algebras in \cite{Schmidtdistr}. We are thus concerned with extending their results to the non-discretely valued case.

Let us briefly recall the proof strategy of \cite[Theorem 3.3]{AWannals} -- the proof of \cite{Schmidtdistr} follows a similar philosophy. 

Suppose that $R$ is a discrete valuation ring and let $\pi\in R$ denote a uniformizer. In order to show Auslander regularity of the Banach algebra in question, consider the unit ball $\A$ (i.e. $\A=R\langle x, \pi^n\partial\rangle$ or $\A=\h{U_R(\pi^n \mathfrak{g}_R)}$, respectively) and pass to the residue field. The resulting $k$-algebra $\A/\pi \A$ carries a natural Zariskian filtration such that $\mathrm{gr}(\A/\pi\A)$ is a polynomial ring over $k$, so that it can be shown via standard techniques that $\A/\pi \A$ is Auslander regular. Thus $\A\otimes_RK$ is also Auslander regular, by invoking the same filtration arguments once more, this time relative to the $\pi$-adic filtration on $\A\otimes_R K$. This explains the terminology of `doubly filtered $K$-algebras' in \cite{AWannals}.

We would like to employ similar ideas in the non-discrete case: elements of our Banach algebras are power series expressions with some convergence condition, but if we work over $R/\pi$ for some $\pi\in \mathfrak{m}\setminus\{0\}$ instead, these become polynomial expressions, where the toolbox of Zariskian filtrations is available to us. Once we have `Auslander regularity mod $\pi$', we would like to lift this back to the setting over $K$.

There is however no hope of this approach working directly: note that $R/\pi$ is itself not even Noetherian. We could try and work over the residue field $k$ instead, but now there is no straightforward way to lift this back to characteristic zero.

It is here that the framework of `almost mathematics' (cf. \cite{GRalmost}) comes to the rescue: considering $R$-modules only `up to modules killed by $\mathfrak{m}$' is routinely used to build a bridge between algebraic information over $R/\pi$ and $K$, when $K$ is not discretely valued. While the algebras over $R$ and $R/\pi$ which we consider are not Noetherian, they are in fact almost Noetherian, and it turns out that we can develop what amounts to an `almost analogue' of the proof in the discretely valued case.

To this end, we introduce the notion of an almost Auslander regular algebra and deduce Theorem \ref{MainThm1} as a combination of the following results:

\begin{thm}
	\label{aAregmodpi}
	Let $\A$ be a $\pi$-adically complete, $\pi$-torsionfree $R$-algebra such that $\A/\pi \A$ is an almost Auslander regular $R/\pi$-algebra. Then $A:=\A\otimes_RK$ is Auslander regular.
\end{thm}

\begin{thm}
	\label{polyaAreg}
	The polynomial ring $R/\pi[x_1, \hdots, x_m]$ is an almost Auslander regular $R/\pi$-algebra for all $m\geq 0$.
\end{thm}

We remark that the Auslander regularity of Tate-Weyl algebras (example (ii)) can be used to establish the Auslander regularity of more general Banach algebras of $p$-adic differential operators on smooth affinoid spaces using suitable Kashiwara equivalences. We will spell this out in a follow-up paper.

This paper is organized as follows:

In section 2, we give a brief overview of the framework of almost mathematics, and develop the necessary tools to use filtered and graded methods in an almost context. 

In section 3, we introduce the notion of an almost Auslander regular $R/\pi$-algebra and prove Theorem \ref{aAregmodpi}.

In section 4, we show that if $\A$ is an almost Auslander regular $R/\pi$-algebra, then $\A[x]$ is almost Auslander regular, provided one extra condition is satisfied: as we do not know of an almost version of Hilbert's Basis Theorem (cf. \cite[Warning 2.7.9]{BZ}), we need to \emph{assume} that $\A[x, y]$ is almost Noetherian. The proof then follows the same strategy as in the classical case in \cite[chapter III, section 2.3]{LVO}: we show that it suffices to consider graded modules and then use various reduction steps so that we only need to consider $\A[x]$-modules which are either already almost finitely generated over $\A$ with $x$ acting nilpotently, or of the form $\A[x]\otimes_\A\M$ for some almost finitely generated $\A$-module $\M$. This allows us to deduce Theorem \ref{polyaAreg} by an inductive argument.

In section 5, we return to the three examples mentioned at the beginning of the introduction and prove their Auslander regularity. To finish the proof of Theorem \ref{MainThm1}, we note that sections 3 and 4 already provide natural bounds for the global dimension. This bound turns out to be sharp in the case of a completed enveloping algebra, so that $\mathrm{gl.dim.}U_n=\mathrm{dim}_K \mathfrak{g}$, which also allows us to determine the global dimension in the case of a distribution algebra. In the appendix, we prove Bernstein's inequality for Tate-Weyl algebras by adapting the proof from \cite{AWannals} to the non-discrete setting, which verifies that $\mathrm{gl.dim.}\h{A_{m, n}}=m$.

\section{Background on almost mathematics}
\subsection{Basic definitions}
We recall some results on almost mathematics. The main references for us are \cite{GRalmost}, \cite{GRfoundations} and \cite{BZ}. Most of the results in this subsection are well-known, even if our presentation differs slightly and e.g. allows for non-commutative rings.

Suppose that $K$ is densely valued, so that $R$ is a (non-Noetherian) non-discrete valuation ring of rank one, and its maximal ideal $\mathfrak{m}$ satisfies $\mathfrak{m}^2=\mathfrak{m}$. Let $\pi\in \mathfrak{m}$ be non-zero.

Let $\A$ be an $R$-algebra. While the literature usually assumes $\A$ to be commutative, we also allow for non-commutative $\A$ in our discussion.

\begin{defn}
	A left $\A$-module $\M$ is called \textbf{almost zero} if $\mathfrak{m}\M=0$. 
\end{defn} 

\begin{defn}
	A morphism $f: \M\to \N$ of left $\A$-modules is called \textbf{almost injective} (resp., \textbf{almost surjective}) if $\mathrm{ker}f$ is almost zero (resp. $\mathrm{coker}f$ is almost zero).
	We say that $f$ is an \textbf{almost isomorphism} if it is almost injective and almost surjective.
\end{defn}

We use the notation $\M\cong^a \N$ to denote that $\M$ is almost isomorphic to $\N$. In particular $\M\cong^a 0$ means that $\M$ is almost zero.

The following is an immediate consequence of the fact that $\mathfrak{m}^2=\mathfrak{m}$.
\begin{lem}
	\leavevmode
	\label{almostsurjcomp}
	\begin{enumerate}[(i)]
		\item The composition of two almost surjective morphisms is almost surjective.
		\item The composition of two almost injective morphisms is almost injective.
		\item If $\M_1\to \M_2\to \M_3$ is an exact sequence of $\A$-modules with $\M_1, \M_3$ almost zero, then $\M_2$ is almost zero.
	\end{enumerate}
\end{lem}

In particular, the subcategory of almost zero left $\A$-modules forms a Serre subcategory of the category $\mathrm{Mod}(\A)$ of all left $\A$-modules. Almost mathematics is concerned with the study of the resulting quotient category, usually denoted $\mathrm{Mod}^a(\A)$. The natural functor $\mathrm{Mod}(\A)\to \mathrm{Mod}^a(\A)$ is denoted by $\M\mapsto \M^a$.

Note that a morphism $f: \M\to \N$ is almost injective (almost surjective, an almost isomorphism) if and only if $f^a: \M^a\to \N^a$ is a monomorphism (an epimorphism, an isomorphism) in $\mathrm{Mod}^a(\A)$.

\begin{defn}
	A left $\A$-module $\M$ is \textbf{almost finitely generated} if for any $\epsilon\in \mathfrak{m}$, there exists a positive integer $n_\epsilon$ and an $\A$-module morphism
	\begin{equation*}
		f_\epsilon: \A^{n_\epsilon}\to \M
	\end{equation*}
	such that $\mathrm{coker}(f_\epsilon)$ is annihilated by $\epsilon$.
\end{defn}

\begin{defn}
	We say that $\A$ is \textbf{almost left Noetherian} if any submodule of a finitely generated left $\A$-module $\M$ is almost finitely generated.
\end{defn}

\begin{lem}
	The following are equivalent:
	\begin{enumerate}[(i)]
		\item $\A$ is almost left Noetherian.
		\item Any left ideal of $\A$ is almost finitely generated.
		\item Any submodule of an almost finitely generated left $\A$-module $\M$ is almost finitely generated.
	\end{enumerate}
\end{lem}
\begin{proof}
	It is obvious that (i) implies (ii) and that (iii) implies (i). It follows from \cite[Corollary 2.7.6]{BZ} that (ii) implies (iii). 
\end{proof}

We define almost right Noetherian analogously, and use the term `almost Noetherian' for $R$-algebras which are almost left Noetherian and almost right Noetherian.

\begin{prop}
	For any $m\geq 0$, the $R$-algebra $R\langle x_1, \hdots, x_m\rangle$ is almost Noetherian.
\end{prop}
\begin{proof}
	This is \cite[Satz 5.1]{Kiehl}.
\end{proof}

It is a common theme that `almost' properties over $R/\pi$ can often be lifted to $R/\pi^2$ (this is already hinted at in Lemma \ref{almostsurjcomp}.(iii)), and then inductively to $R/\pi^n$ for any $n$. Under the assumption of $\pi$-adic completeness, this will often allow us to lift almost properties from $R/\pi$ to $R$.

We first record the following general lemma:
\begin{lem}
	\label{afgcompletion}
	Let $\A$ be an $R$-algebra, and let $\M$ be a $\pi$-adically separated left $\A$-module. Fix $\epsilon\in \mathfrak{m}$ with $|\epsilon|>|\pi|$, and let $\N\subseteq \M$ be a finitely generated $\A$-submodule such that
	\begin{equation*}
		\epsilon \M\subseteq \N+\pi\M.
	\end{equation*}
	\begin{enumerate}[(i)]
		\item We have 
		\begin{equation*}
			\epsilon\cdot \M\subseteq \h{\A}\cdot \N\subseteq \h{\M}.
		\end{equation*}
		In particular, if $\A$ is $\pi$-adically complete, then a $\pi$-adically separated left $\A$-module $\M$ is almost finitely generated as soon as $\M/\pi \M$ is almost finitely generated.
		\item If $\M$ is $\pi$-torsionfree, then the map $\h{\N}\to \h{\M}$ is injective, and $\h{\N}\otimes_RK\to \h{\M}\otimes_RK$ is an isomorphism.
	\end{enumerate}
\end{lem}
\begin{proof}
	\begin{enumerate}[(i)]
		\item Let $x_1, \hdots, x_r\in \M$ be generators of $\N$, such that
		\begin{equation*}
			\epsilon\cdot\M\subseteq \sum \A\cdot x_i+\pi \M.
		\end{equation*}
		If $m_0\in \M$, then we can find $a_{0, 1}, \hdots, a_{0,r}\in \A$ and $m_1\in \M$ such that
		\begin{equation*}
			\epsilon m_0=\sum a_{0,i}x_i+\pi \epsilon^{-1}(\epsilon m_1).
		\end{equation*}
		Applying the same argument to $m_1$, we have $a_{1, 1}, \hdots, a_{1, r}\in \A$ and $m_2\in \M$ such that
		\begin{equation*}
			\epsilon m_0=\sum (a_{0,i}+\pi \epsilon^{-1}a_{1, i})x_i+(\pi \epsilon^{-1})^2 (\epsilon m_2).
		\end{equation*}
		Iterating this process, we note that for each $i$, the partial sums $\sum_{j=0}^n (\pi \epsilon^{-1})^ja_{j, i}$ converge to an element $\lambda_i\in \h{\A}$. Since $\M$ is $\pi$-adically separated, it embeds into $\h{\M}$, and now $\epsilon m_0-\sum \lambda_ix_i\in \cap \pi^n\h{\M}=\{0\}$. Thus $\epsilon m_0\in \sum \h{\A}x_i\subseteq \h{\M}$, as required.
		In particular, if $\A$ is $\pi$-adically complete and $\M$ is a $\pi$-adically separated left $\A$-module such that $\M/\pi\M$ is almost finitely generated, then we can choose for any $\epsilon\in \mathfrak{m}$ with $|\epsilon|>|\pi|$ finitely many elements $x_1, \hdots, x_r\in \M$ such that
		\begin{equation*}
			\epsilon \M\subseteq \sum_i \A\cdot x_i+\pi \M.
		\end{equation*}
		The above then shows that $\epsilon\M\subseteq \sum_i\A\cdot x_i$, and hence $\M$ is almost finitely generated.
		\item Let $x\in \h{\N}$ be an element in the kernel of the map $\h{\N}\to \h{\M}$. Thus there exists $x_n\in \N$ with $x-x_n\in \pi^n\h{\N}$ and $x_n\in \N\cap \pi^n\M$. Write $x_n=\pi^ny_n$ with $y_n\in \M$. By assumption, we have
		\begin{equation*}
			\epsilon\cdot y_n\in \N+\pi \M\subseteq \h{\N}+\pi\M\subseteq \h{\N}\otimes_RK.
		\end{equation*} 
		Repeating the same argument as in (i), using that $\h{\N}$ is an $\h{\A}$-module and $\M$ is $\pi$-adically separated, it follows that $\epsilon y_n\in \h{\N}$. But then $\epsilon x_n\in \pi^n\h{\N}$, so $\epsilon x\in\cap \pi^n\h{\N}=\{0\}$. Since $\N\subseteq \M$ is $\pi$-torsionfree, so is $\h{\N}$, and therefore $x=0$. Hence $\h{\N}\to \h{\M}$ is injective. Finally, it then follows by the same strategy as in (i) that $\epsilon \h{\M}\subseteq \h{\N}$, since $\h{\M}$ is $\pi$-adically separated. Thus, the cokernel of the map $\h{\N}\to \h{\M}$ is annihilated by $\epsilon$, and the map becomes an isomorphism after tensoring with $K$.\qedhere
	\end{enumerate}
\end{proof}

\begin{cor}
	\label{bdedptormodp}
	Let $\A$ be an $R$-algebra such that $\A/\pi\A$ is almost Noetherian. Then any almost finitely generated, $\pi$-adically separated left $\A$-module has bounded $\pi$-torsion, i.e. if $\M$ is almost finitely generated and $\pi$-adically separated, then there exists an integer $s$ such that
	\begin{equation*}
		\pi^nx=0 \text{ for some $n$} \Rightarrow \pi^sx=0
	\end{equation*}
	for all $x\in \M$.
\end{cor}
\begin{proof}
	Let $T=\{x\in \M: \pi^nx=0 \text{ for some $n$}\}$. Since $M/T$ is $\pi$-torsionfree, we have $T/\pi T\subseteq \M/\pi\M$, so that $T/\pi T$ is almost finitely generated by assumption.
	
	If $\epsilon\in \mathfrak{m}$ with $|\epsilon|>|\pi|$, this yields a finitely generated $\A$-submodule $T'\subseteq T$ such that
	\begin{equation*}
		\epsilon T\subseteq T'+\pi T.
	\end{equation*}
	Applying Lemma \ref{afgcompletion}.(i), we have $\epsilon\cdot T\subseteq \h{\A}\cdot T'$. Since $T'\subseteq T$ is finitely generated over $\A$, it is annihilated by some $\pi^{s'}$, and thus $T$ is annihilated by $\epsilon\cdot \pi^{s'}$.
\end{proof}

The following result now provides a first illustration of the principle `almost properties over $R/\pi$ yield properties over $K$'.

\begin{prop}
	\label{liftNoeth}
	Let $\A$ be a $\pi$-adically complete $R$-algebra, and suppose that $\A/\pi\A$ is almost left Noetherian. Then $\A$ is almost left Noetherian, and $A:=\A\otimes_RK$ is left Noetherian.
\end{prop}

The proof rests on the following two lemmas:

\begin{lem}
	Let $\A$ be as in the statement of Proposition \ref{liftNoeth}.
	
	Let $I\subseteq \A$ be a left ideal. Suppose that there exists an almost finitely generated left ideal $J$ and some $s\geq 0$ such that $\pi^sJ\subseteq I\subseteq J$. Then $I$ is almost finitely generated.
\end{lem}

\begin{proof}
	We argue by induction on $s$, the case $s=0$ being trivial.
	
	For $s\geq 1$, let $J_0=\{x\in J: \pi x\in I\}$, a left ideal of $\A$. We have
	\begin{equation*}
		I\subseteq J_0\subseteq J,
	\end{equation*}
	and $\pi^{s-1}J\subseteq J_0$ by construction. Thus $J_0$ is almost finitely generated by the induction hypothesis.
	
	We now have a short exact sequence of $\A$-modules
	\begin{equation*}
		0\to \pi J_0\to I\to I/\pi J_0\to 0.
	\end{equation*}
	Since $I/\pi J_0$ is a submodule of the almost finitely generated $\A/\pi \A$-module $J_0/\pi J_0$, it is almost finitely generated over $\A/\pi \A$, and hence almost finitely generated over $\A$. As $\pi J_0$ is also almost finitely generated, the result follows by Lemma \ref{almostsurjcomp}.(iii). 
\end{proof}

\begin{lem}
	\label{satisafg}
	Let $\A$ be as in the statement of Proposition \ref{liftNoeth}.
	
	Let $I\subseteq \A$ be a left ideal, and let $J=\{x\in \A: \pi^nx\in I \ \text{for some } n\}$ be its saturation.
	\begin{enumerate}[(i)]
		\item $J$ is almost finitely generated.
		\item There exists some $s\geq 0$ such that $\pi^sJ\subseteq I$.
	\end{enumerate}
\end{lem}

\begin{proof}
	\begin{enumerate}[(i)]
		\item Since $J\cap \pi\A=\pi J$, it follows that $J/\pi J\subseteq \A/\pi\A$ is a left ideal and hence almost finitely generated by assumption. The statement thus follows directly from Lemma \ref{afgcompletion}.(i).
		\item Let $\epsilon\in \mathfrak{m}$, $|\epsilon|>|\pi|$. By (i), there exist $x_1, \hdots, x_r\in J$ such that $\epsilon J\subseteq \sum \A x_i$. By definition of $J$, there exists $s_i$ such that $\pi^{s_i}x_i\in I$, so in particular $\pi^t(\sum \A x_i)\subseteq I$ for any $t\geq \mathrm{max}\{s_i\}$. Taking $s=\mathrm{max}\{s_i\}+1$ proves the result.\qedhere 
	\end{enumerate}
\end{proof}

\begin{proof}[Proof of Proposition \ref{liftNoeth}]
	It follows immediately from the two lemmas above that $\A$ is almost Noetherian. 
	
	Thus, the image of $\A$ in $A=\A\otimes_RK$ is likewise almost Noetherian, and we suppose without loss of generality that $\A$ is $\pi$-torsionfree, so that $\A\subseteq A$.
	
	If $I$ is a left ideal in $A$, let $\mathcal{I}=\A\cap I$, which is almost finitely generated by the above. For any $\epsilon \in \mathfrak{m}$, there exists a finitely generated ideal $J\subseteq \mathcal{I}$ such that $\epsilon \mathcal{I}\subseteq J$. In particular, $I=\mathcal{I}\otimes_RK\cong J\otimes_RK$ is a finitely generated ideal in $A$.
\end{proof}

If $\A$ is a $\pi$-adically complete almost Noetherian $R$-algebra, then we have an analogue of the usual Artin-Rees theory.

\begin{lem}
	\label{almostAR}
	Let $\A$ be a $\pi$-adically complete and almost Noetherian $R$-algebra. Let $\M$ be a finitely generated left $\A$-module.
	\begin{enumerate}[(i)]
		\item If $\N\subseteq \M$ is a left submodule, then there exists an $n\geq 0$ such that
		\begin{equation*}
			\N\cap \pi^{n+1}\M=\pi(\N\cap \pi^n\M).
		\end{equation*}
		\item $\M$ is $\pi$-adically separated.
		\item $\M$ is $\pi$-adically complete.
	\end{enumerate}
\end{lem}

\begin{proof}
	\begin{enumerate}[(i)]
		\item Let $\N^{\mathrm{sat}}=\{x\in \M: \pi^sx\in \N \ \text{for some }\ s\}$ be the saturation of $\N$. By the same argument as in Lemma \ref{satisafg}, there exists some $n$ such that $\pi^n\N^{\mathrm{sat}}\subseteq \N$, and hence
		\begin{equation*}
			\pi^{n+1}\M\cap \N=\pi^{n+1}\N^{\mathrm{sat}}\cap \N= \pi(\pi^n\M\cap \N),
		\end{equation*}
		as required.
		\item If $x\in \cap \pi^n\M$, let $\N=\A\cdot x$. By the above, $\N=\pi \N$, but then there exists $a\in \A$ such that $(1-\pi a)\cdot x=0$. Since $1-\pi a$ is a unit in $\A$, it follows that $x=0$.
		\item Suppose first that $\M$ is $\pi$-torsionfree, and consider the short exact sequence
		\begin{equation*}
			0\to \N\to\A^r\to \M\to 0.
		\end{equation*}
		By assumption, this yields the short exact sequence
		\begin{equation*}
			0\to \varprojlim \frac{\N}{\N\cap \pi^n \A^r}\to \A^r\to \h{\M}\to 0,
		\end{equation*}
		since $\A$ is $\pi$-adically complete. We remark that the Artin-Rees property in (i) ensures that the first term is actually isomorphic to $\h{\N}$, but the above is already to enough to observe that the natural morphism $\M\to \h{\M}$ is surjective, as the surjection $\A^r\to \h{\M}$ factors through $\M$. Since $\M$ is $\pi$-adically separated by point (ii), this shows that $\M\to \h{\M}$ is an isomorphism, i.e. $\M$ is $\pi$-adically complete if it is $\pi$-torsionfree.
		
		In general, let $T$ denote the submodule of torsion elements in $\M$. Note that this is the saturation of $0$ in $\M$, and hence there exists $s$ such that $\pi^sT=0$, by the same argument as before. In particular, $T$ is $\pi$-adically complete. 
		
		But then $\M$ is the extension of the two $\pi$-adically complete $\A$-modules $T$ and $\M/T$ (which is $\pi$-adically complete since it is $\pi$-torsionfree and finitely generated). The Artin-Rees property and $\M/T$ being $\pi$-torsionfree ensures once more that we have a commutative diagram
		\begin{equation*}
			\begin{xy}
				\xymatrix{0\ar[r]& T\ar[r]\ar[d]& \M\ar[r]\ar[d]& \M/T\ar[r]\ar[d]&0\\
					0\ar[r]&\h{T}\ar[r]&\h{\M}\ar[r]&\h{\M/T}\ar[r]&0}
			\end{xy}
		\end{equation*}
		where the outer vertical arrows are isomorphisms. Hence the morphism in the middle is also an isomorphism and $\M$ is $\pi$-adically complete.
	\end{enumerate}
\end{proof}

\begin{cor}
	Let $\A$ be a $\pi$-adically complete, almost Noetherian $R$-algebra. Any almost finitely generated $\A$-module has bounded $\pi$-torsion, i.e. if $\M$ is almost finitely generated, then there exists an integer $s$ such that
	\begin{equation*}
		\pi^nx=0 \text{ for some $n$} \Rightarrow \pi^sx=0
	\end{equation*}
	for all $x\in \M$.
\end{cor}
\begin{proof}
	It is straightforward to reduce the claim to the case where $\M$ is finitely generated. But then $\M$ is $\pi$-adically separated by Lemma \ref{almostAR}.(ii), so we are done by Corollary \ref{bdedptormodp}.
\end{proof}

\begin{cor}
	\label{ptfisfp}
	Let $\A$ be a $\pi$-adically complete, almost Noetherian $R$-algebra. A finitely generated left $\A$-module which is $\pi$-torsionfree is finitely presented.
	
	In particular, a $\pi$-adically saturated submodule of a finitely generated left $\A$-module is itself finitely generated.
\end{cor}
\begin{proof}
	Let $\M$ be a $\pi$-torsionfree, finitely generated left $\A$-module, and consider the short exact sequence
	\begin{equation*}
		0\to \N\to \A^r\to \M\to 0.
	\end{equation*} 
	Since $\M$ is $\pi$-torsionfree, it is flat over $R$, so
	\begin{equation*}
		0\to \N/\pi\N\to (\A/\pi\A)^r\to \M/\pi\M\to 0
	\end{equation*}
	is also exact. Now $\M/\pi\M$ is finitely generated over $\A/\pi\A$ and flat over $R/\pi$, so it is finitely presented by \cite[Theorem 3.4.6]{RG}. 
	
	It thus follows from \cite[tag 0519]{stacksproj} that $\N/\pi\N$ is a finitely generated left $\A/\pi$-module. Let $x_1, \hdots, x_s\in \N$ generate $\N/\pi\N$. Since $\M$ is flat over $R$, we also obtain the exact sequence
	\begin{equation*}
		0\to \h{\N}\to \h{\A}^r\to \h{\M}\to 0,
	\end{equation*} 
	and since both $\A$ and $\M$ are $\pi$-adically complete by Lemma \ref{almostAR}.(iii), we can deduce that $\N$ is $\pi$-adically complete. This now allows us to use a similar argument to Lemma \ref{afgcompletion}: if $n_0\in \N$, there exist $a_{0, 1}, \hdots, a_{0, s}\in \A$ such that
	\begin{equation*}
		n_0=\sum a_{0, i}x_i+\pi n_1
	\end{equation*} 
	for some $n_1\in \N$. Iterating this process, let $\lambda_i\in \A$ denote the limit of the partial sums $\sum_{j=0}^n \pi^j a_{j, i}$ to deduce that $n_0=\sum \lambda_i x_i\in \sum \A\cdot x_i$, so that $\N$ is finitely generated.
	
	The last claim follows immediately: If $\M$ is a finitely generated left $\A$-module and $\N\subseteq \M$ is $\pi$-adically saturated, then $\M/\N$ is finitely generated and $\pi$-torsionfree, and thus finitely presented by the above. We can thus deduce from \cite[tag 0519]{stacksproj} that $\N$ is finitely generated.
\end{proof}

If $\A$ is an almost Noetherian $R$-algebra, then the same argument as in \cite[Proposition 4.3.4]{Hausdf} ensures that the morphism $\A\to \h{\A}$ is flat. We conclude this subsection with a partial generalization of this result, which turns out to be particularly useful for proving flatness statements when working with Fr\'echet--Stein algebras.

\begin{prop}
	Let $\A$ be an almost Noetherian, $\pi$-torsionfree $R$-algebra, and let $\B\subseteq \A\otimes_RK$ be a $\pi$-adically separated $R$-subalgebra containing $\A$. Suppose that $\B/\pi\B$ is almost Noetherian. Then the natural morphism 
	\begin{equation*}
		\A\otimes_RK=\B\otimes_RK\to \h{\B}\otimes_RK
	\end{equation*}
	is flat on both sides.
\end{prop}
\begin{proof}
	We first prove that if $J\subseteq \B$ is a finitely generated left ideal, then the natural morphism $\h{\B}\otimes_{\B}J\to \h{J}$ is an almost isomorphism.
	
	In fact, consider a short exact sequence
	\begin{equation*}
		0\to \N\to \B^r\to J\to 0,
	\end{equation*}
	giving rise to a commutative diagram
	\begin{equation*}
		\begin{xy}
			\xymatrix{& \h{\B}\otimes_{\B}\N\ar[r]\ar[d]&\h{\B}^r\ar[r]\ar[d]&\h{\B}\otimes_{\B}J\ar[r]\ar[d]&0\\
				0\ar[r]&\h{\N}\ar[r]&\h{\B}^r\ar[r]&\h{J}\ar[r]&0}
		\end{xy}
	\end{equation*}
	with exact rows, since $J$ is $\pi$-torsionfree. In particular, the morphism $\h{\B}\otimes_{\B}J\to \h{J}$ is surjective.
	
	By Proposition \ref{liftNoeth}, $\h{\B}$ is almost Noetherian, so both $\h{J}$ and $\h{\N}$ are almost finitely generated. Since $\N$ is $\pi$-adically separated, we can now use the same strategy as in Lemma \ref{afgcompletion} to deduce that the map $\h{\B}\otimes_{\B}\N\to \h{\N}$ is almost surjective, and thus $\h{\B}\otimes_{\B}J\to \h{J}$ is almost injective, as required.
	
	In particular, tensoring with $K$ yields an isomorphism $(\h{\B}\otimes_{\B}J)\otimes_RK\to\h{J}\otimes_RK$ for any finitely generated left ideal $J\subseteq \B$.
	
	Now write $A=\A\otimes_RK=\B\otimes_RK$, $\h{B}=\h{\B}\otimes_RK$.
	
	Let $I\subseteq A$ be a left ideal. We need to show that the morphism $\h{B}\otimes_AI\to \h{B}$ is injective.
	
	Since $\B$ is $\pi$-adically separated, it embeds into $\h{\B}$, so the map $A\to \h{B}$ is injective. 
	
	Let $\I=I\cap \B$, a $\pi$-adically saturated left ideal of $\B$. In particular, writing $S=\B/\I$, the short exact sequence
	\begin{equation*}
		0\to \I\to \B\to S\to 0
	\end{equation*}
	remains exact after $\pi$-adic completion, since $S$ is $\pi$-torsionfree. Thus, $\h{\I}$ is a left ideal of $\h{\B}$, and it suffices to show that $\h{B}\otimes_AI=\h{B}\otimes_{\B} \I\to \h{\I}\otimes_RK$ is injective. Since $\h{\B}$ is almost Noetherian, we know that $\h{\I}$ is almost finitely generated. In particular, $\I/\pi \I\cong \h{\I}/\pi \h{\I}$ is an almost finitely generated $\B$-module.
	
	Fix $\epsilon\in \mathfrak{m}$ with $|\epsilon|>|\pi|$. Let $J\subseteq \I$ be a finitely generated left ideal of $\B$ such that
	\begin{equation*}
		\epsilon \I\subseteq J+\pi \I.
	\end{equation*}
	Since $\I\cap \A$ is a left ideal of the almost Noetherian algebra $\A$, we can suppose without loss of generality that $\epsilon (\I\cap \A)\subseteq J$. In particular, for any $x\in \I\subseteq \A\otimes_RK$, there exists some $n$ such that $\pi^nx\in J$, and thus $J\otimes_RK\cong\I\otimes_RK=I$.
	
	Since $\I$ is $\pi$-torsionfree and $\pi$-adically separated, we can apply Lemma \ref{afgcompletion}.(ii) to deduce that $\h{J}\to \h{\I}$ is injective, and $\h{J}\otimes_RK\to \h{\I}\otimes_RK$ is an isomorphism.
	
	We thus have the following commutative diagram
	\begin{equation*}
		\begin{xy}
			\xymatrix{\h{B}\otimes_{\B}J\ar[r]\ar[d]&\h{B}\otimes_{\B}\I\ar[d]\\
				\h{J}\otimes_RK\ar[r]&\h{\I}\otimes_RK}
		\end{xy}
	\end{equation*}
	with the left vertical map an isomorphism (by the first step in our proof) and the bottom horizontal map an isomorphism.
	
	We have already remarked above that $J\otimes_RK\cong \I\otimes_R K$, so the top horizontal arrow is in fact an isomorphism, as it can be written as the natural morphism $\h{B}\otimes_A(J\otimes_RK)\to \h{B}\otimes_A I$. It thus follows that $\h{B}\otimes_{\B} \I\to \h{\I}\otimes_RK$ is an isomorphism, which finishes the proof.
\end{proof}

\subsection{Good filtrations and almost mathematics}
In this subsection, we develop the almost analogue of the theory of good filtrations and related graded methods. We refer to \cite{LVO} for an overview of the classical results.

Let $\A$ be an almost Noetherian $R$-algebra, equipped with an exhaustive filtration
\begin{equation*}
	0=\mathrm{Fil}_{-1}\A\subseteq \mathrm{Fil}_0\A\subseteq \mathrm{Fil}_1\A\subseteq \hdots \subseteq \A,
\end{equation*}
endowing $\A$ with the structure of a filtered $R$-algebra, i.e. each $\mathrm{Fil}_i\A$ is an $R$-submodule, and we have $\mathrm{Fil}_i\A\cdot \mathrm{Fil}_j\A\subseteq \mathrm{Fil}_{i+j}\A$ for all $i, j$. We also suppose that $1\in \mathrm{Fil}_0\A$.

\begin{defn}
	A filtered $R$-algebra $\A$ as above is called \textbf{almost Zariskian} if the associated Rees ring
	\begin{equation*}
		\mathrm{Rees}\A=\oplus_i \mathrm{Fil}_i\A\cdot t^i\subseteq \A[t]
	\end{equation*}
	is almost Noetherian.
\end{defn}
Throughout this subsection, $\A$ will denote an almost Zariskian $R$-algebra. Note that both $\A\cong \frac{\mathrm{Rees} \A}{(t-1)\mathrm{Rees} \A}$ and the associated graded ring
\begin{equation*}
	\mathrm{gr}\A=\underset{i}{\bigoplus} \frac{\mathrm{Fil}_i\A}{\mathrm{Fil}_{i-1}\A}\cong \frac{\mathrm{Rees} \A}{t\cdot \mathrm{Rees} \A}
\end{equation*}
are quotients of $\mathrm{Rees}\A$ and hence almost Noetherian.

An $\A$-module filtration on an $\A$-module $\M$ is a filtration
\begin{equation*}
	\hdots\subseteq \mathrm{Fil}_i\M\subseteq \mathrm{Fil}_{i+1}\M\subseteq \hdots \subseteq \M
\end{equation*}
by $R$-submodules satisfying $\mathrm{Fil}_i\A\cdot \mathrm{Fil}_j\M\subseteq \mathrm{Fil}_{i+j}\M$ for all $i, j$. All our filtrations will be left bounded in the sense that $\mathrm{Fil}_i\M=0$ for $i<<0$, and exhaustive, meaning $\cup_i \mathrm{Fil}_i\M=\M$.

Analogously to the above, we can form the Rees module $\mathrm{Rees}\M$ and the graded module $\mathrm{gr}\M$, which are left modules over $\mathrm{Rees}\A$ and $\mathrm{gr}\A$, respectively.

\begin{defn}
	Let $\M$ be an $\A$-module. An $\A$-module filtration on $\M$ is called an \textbf{almost good filtration} if $\mathrm{Rees}\M$ is an almost finitely generated $\mathrm{Rees}\A$-module.
\end{defn}

In the classical setting, it is more common to require $\mathrm{gr}\M$ to be finitely generated, which however turns out to be equivalent to $\mathrm{Rees}\M$ being finitely generated (see \cite[chapter I, Lemma 5.3]{LVO}).

For our purposes, the Rees construction will be crucial, making the definition as stated more convenient. It is clear that if $\M$ carries an almost good filtration, then $\M$ is almost finitely generated over $\A$ and $\mathrm{gr}\M\cong \mathrm{Rees}\M/t\mathrm{Rees}\M$ is almost finitely generated over $\mathrm{gr}\A$. It is not clear (and quite possibly false) that a filtration with $\mathrm{gr}\M$ being almost finitely generated is always almost good as defined above.

Note that if $\M$ is an $\A$-module endowed with an almost good filtration and $\N\subseteq \M$ is a submodule, then the induced subspace filtration on $\N$ is an almost good filtration, as $\mathrm{Rees}\N\subseteq \mathrm{Rees}\M$ and $\mathrm{Rees}\A$ is almost Noetherian. The corresponding statement for quotient filtrations is trivially true.

\begin{defn}
	Let $\M$, $\N$ be filtered $\A$-modules. Let $\epsilon\in \mathfrak{m}$.
	\begin{enumerate}[(i)]
		\item An $\A$-module morphism $f: \M\to \N$ is called a \textbf{morphism of filtered $\A$-modules} if $f(\mathrm{Fil}_i\M)\subseteq \mathrm{Fil}_i\N$ for each $i$.
		\item A morphism $f: \M\to \N$ of filtered $\A$-modules is called \textbf{strict} if for each $i$,
		\begin{equation*}
			\mathrm{Im}f\cap \mathrm{Fil}_i\N=f(\mathrm{Fil}_i\M).
		\end{equation*}
		\item A morphism $f: \M\to \N$ of filtered $\A$-modules is called \textbf{$\epsilon$-strict} if for each $i$, 
		\begin{equation*}
			\epsilon\cdot (\mathrm{Im}f\cap \mathrm{Fil}_i\N)\subseteq f(\mathrm{Fil}_i\M).
		\end{equation*}
	\end{enumerate}
\end{defn}

Note that a morphism $f:\M\to\N$ of filtered $\A$-modules induces a morphism $\mathrm{gr}f: \mathrm{gr}\M\to \mathrm{gr}\N$. If we endow $\mathrm{Im}f$ with the subspace filtration, then $\mathrm{gr}(\mathrm{Im}f)$ is a submodule of $\mathrm{gr}\N$ containing $\mathrm{Im}(\mathrm{gr}f)$. If $f$ is $\epsilon$-strict, then
\begin{equation*}
	\epsilon\cdot \frac{\mathrm{gr}(\mathrm{Im}f)}{\mathrm{Im}(\mathrm{gr}f)}=0.
\end{equation*}
Similarly, endowing $\mathrm{ker}f$ with the subspace filtration, $\mathrm{gr}(\mathrm{ker}f)$ is a submodule of $\mathrm{gr}\M$ which is contained in $\mathrm{ker}(\mathrm{gr}f)$. If $m\in \mathrm{Fil}_n\M\setminus \mathrm{Fil}_{n-1}\M$ such that $f(m)\in \mathrm{Fil}_{n-1}\M$ (so that its principal symbol $\sigma(m)\in \mathrm{gr}_n\M$ is in the kernel of $\mathrm{gr}f$), the condition of $\epsilon$-strictness also implies that there exists $m'\in \mathrm{Fil}_{n-1}\M$ such that $\epsilon m-m'\in \mathrm{ker}f$. Therefore, we also have
\begin{equation*}
	\epsilon\cdot \frac{\mathrm{ker}(\mathrm{gr}f)}{\mathrm{gr}(\mathrm{ker}f)}=0
\end{equation*}
if $f$ is $\epsilon$-strict.

\begin{defn}
	Let $\epsilon\in \mathfrak{m}$. A sequence of $\A$-modules
	\begin{equation*}
		\begin{xy}
			\xymatrix{\L\ar[r]^f&\M\ar[r]^g& \N
			}
		\end{xy}
	\end{equation*}
	is called \textbf{$\epsilon$-exact} (at $\M$) if $\epsilon \cdot\mathrm{ker}g\subseteq \mathrm{Im}f\subseteq \mathrm{ker}g$.
\end{defn}

\begin{prop}
	\label{Reesapproximation}
	Let $\epsilon\in \mathfrak{m}$. Let $f: \L\to \M$, $g: \M\to \N$ be $\epsilon$-strict morphisms of filtered $\A$-modules such that
	\begin{equation*}
		\begin{xy}
			\xymatrix{\L\ar[r]^f&\M\ar[r]^g&\N
			}
		\end{xy}
	\end{equation*}
	is $\epsilon$-exact at $\M$. Then the resulting sequence
	\begin{equation*}
		\mathrm{Rees}\L\to \mathrm{Rees}\M\to \mathrm{Rees}\N
	\end{equation*}
	is $\epsilon^2$-exact at $\mathrm{Rees}\M$, and the sequence
	\begin{equation*}
		\mathrm{gr}\L\to \mathrm{gr}\M\to \mathrm{gr}\N
	\end{equation*}
	is $\epsilon^3$-exact at $\mathrm{gr}\M$.
\end{prop}
\begin{proof}
	Note that $\mathrm{Rees}g\circ \mathrm{Rees}f=0$ by functoriality. If $m\in \mathrm{Fil}_i\M$ is in the kernel of $g$, then it follows from $\epsilon$-exactness that $\epsilon\cdot m\in \mathrm{Im}f\cap \mathrm{Fil}_i\M$. Hence $\epsilon^2\cdot m\in f(\mathrm{Fil}_i\L)$, since $f$ is $\epsilon$-strict. This shows that the sequence of Rees modules is $\epsilon^2$-exact.
	
	Moreover, $\mathrm{gr}g\circ \mathrm{gr}f=0$ by functoriality. Now
	\begin{equation*}
		\epsilon \mathrm{ker}(\mathrm{gr}g)\subseteq \mathrm{gr}(\mathrm{ker}g)
	\end{equation*} 
	by $\epsilon$-strictness of $g$. Due to $\epsilon$-exactness, we have
	\begin{equation*}
		\epsilon\cdot \mathrm{ker}g\subseteq \mathrm{Im}f
	\end{equation*}
	and hence
	\begin{equation*}
		\epsilon\cdot \mathrm{gr}(\mathrm{ker}g)\subseteq \mathrm{gr}\mathrm{Im}f.
	\end{equation*}
	Lastly, we have
	\begin{equation*}
		\epsilon\cdot \mathrm{gr}(\mathrm{Im}f)\subseteq \mathrm{Im}(\mathrm{gr}f)
	\end{equation*}
	by $\epsilon$-strictness of $f$, and therefore 
	\begin{equation*}
		\epsilon^3\cdot \mathrm{ker}(\mathrm{gr}g)\subseteq \mathrm{Im}(\mathrm{gr}f),
	\end{equation*}
	as required.
\end{proof}

\begin{defn}
	Let $\M^\bullet$ be a chain complex of $\A$-modules. For $\epsilon\in \mathfrak{m}$, we say that a morphism of chain complexes $f: \F^\bullet\to \M^\bullet$ is an \textbf{$\epsilon$-approximation} of $\M^\bullet$ if the cone of $f$ is $\epsilon$-exact. 
\end{defn}
If $f$ is an $\epsilon$-approximation of $\M^\bullet$, then $f$ induces a morphism between the cohomololgy groups
\begin{equation*}
	\mathrm{H}^j(\F^\bullet)\to \mathrm{H}^j(\M^\bullet)
\end{equation*}
whose kernel and cokernel are annihilated by $\epsilon$.

\begin{lem}
	\label{Extapproximation}
	Let $\M$ be an $\A$-module and let $\F^\bullet\to\M$ be an $\epsilon$-approximation for some $\epsilon\in \mathfrak{m}$, where $\F^\bullet$ is a complex of free $\A$-modules concentrated in non-positive degrees. Let
	\begin{equation*}
		\alpha: \mathrm{Ext}^j_\A(\M, \A)\to \mathrm{H}^j(\mathrm{Hom}_\A(\F^\bullet, \A))
	\end{equation*}
	denote the natural morphism. Then $\epsilon^{j+1}\mathrm{ker}\alpha=0$ and $\epsilon^{j+2}\mathrm{coker}\alpha=0$.
\end{lem}
\begin{proof}
	Let $C^\bullet$ be the cone of $\F^\bullet\to \M$, i.e.
	\begin{equation*}
		C^\bullet=(\hdots \to \F^{-1}\to \F^0\to \M\to 0),
	\end{equation*}
	an $\epsilon$-exact complex concentrated in non-positive degrees (with $\F^{-i}$ being placed in degree $-i-1$). Considering the long exact sequence of cohomology corresponding to the distinguished triangle
	\begin{equation*}
		\mathrm{R}\mathrm{Hom}_\A(C^\bullet, \A)\to \mathrm{R}\mathrm{Hom}_\A(\M, \A)\to \mathrm{R}\mathrm{Hom}_\A(\F^\bullet, \A)\to^{+1},
	\end{equation*}
	it suffices to show that
	\begin{equation*}
		\epsilon^{j+1} \mathrm{H}^j(\mathrm{R}\mathrm{Hom}_\A(C^\bullet, \A))=0
	\end{equation*}
	for any $j$.
	
	But this follows straightforwardly from $C^\bullet$ being $\epsilon$-exact: note that the cohomology group $\mathrm{H}^j(\mathrm{R}\mathrm{Hom}_\A(C^\bullet, \A))$ carries a $j+1$-step filtration where each graded piece is a subquotient of $\mathrm{Ext}^p_\A(\mathrm{H}^q(C^\bullet), \A)$ for suitable $p$ and $q$ with $p+q=j$ (either by a spectral sequence argument as in \cite[tag 015J]{stacksproj} or by an inductive argument along the long exact sequence of Ext groups for various truncations of $C^\bullet$).
\end{proof}

For example, if $\M$ is an almost finitely generated $\A$-module, then for each $\epsilon\in \mathfrak{m}$ there exists an $\epsilon$-approximation $f: \F^\bullet\to \M$ with $\F^\bullet$ a chain complex concentrated in non-positive degrees and each $\F^i$ a finite free $\A$-module, due to $\A$ being almost Noetherian (cf. a similar discussion in \cite[subsection 2.6]{BZ}).

Our next aim is to show that if $\M$ is an almost finitely generated $\A$-module endowed with an almost good filtration, we can produce such an $\epsilon$-approximation compatible with filtrations in the sense that it involves only $\epsilon$-strict morphisms.

\begin{defn}
	\leavevmode
	\begin{enumerate}[(i)]	
		\item Let $\M$ be a filtered $\A$-module. For any $n\in \mathbb{Z}$, the \textbf{$n$-shift} of $\M$ is the filtered $\A$-module $\M(n)$ whose underlying $\A$-module is $\M$, with filtration
		\begin{equation*}
			\mathrm{Fil}_i\M(n)=\mathrm{Fil}_{i+n}\M.
		\end{equation*}
		\item A filtered $\A$-module $\F$ is called \textbf{finite filt-free} if it is of the form
		\begin{equation*}
			\F=\oplus_{i=1}^r \A(n_i)
		\end{equation*}
		(or $\F=0$) for some $r$ and some integers $n_1, \hdots, n_r\in \mathbb{Z}$.
	\end{enumerate}
\end{defn}

\begin{prop}
	Let $\M$ be an almost finitely generated $\A$-module equipped with an almost good filtration. For each $\epsilon\in \mathfrak{m}$, there exists an $\epsilon$-strict morphism
	\begin{equation*}
		f: \F\to \M,
	\end{equation*}
	where $\F$ is a finite filt-free $\A$-module and $\epsilon\cdot \mathrm{coker}f=0$.
\end{prop}
\begin{proof}
	By assumption, $\mathrm{Rees}\M$ is almost finitely generated over $\mathrm{Rees}\A$, so we can find $m_1, \hdots, m_r\in \mathrm{Rees}\M$ such that $\epsilon\cdot \mathrm{Rees}\M\subseteq \sum_i \mathrm{Rees}\A\cdot m_i$. Without loss of generality, we can assume that the $m_i$ are homogeneous elements of degree $d_i$, and that $m_i\notin t\cdot \mathrm{Rees}\M$.
	Thus, $m_i$ corresponds to an element $x_i\in \mathrm{Fil}_{d_i}\M\setminus \mathrm{Fil}_{d_i-1}\M$, i.e. $m_i$ is the element $x_i t^{d_i}\in \mathrm{Fil}_{d_i}\M t^{d_i}$. Note in particular that $\epsilon \M\subseteq \sum_i \A\cdot x_i$ by construction.
	
	We thus obtain a morphism of filtered $\A$-modules
	\begin{equation*}
		f:\F:=\oplus_{i=1}^r \A(-d_i)\to \M
	\end{equation*} 
	given by sending $1$ in the $i$th summand to $x_i$.
	
	The natural quotient map
	\begin{equation*}
		\mathrm{Rees}\M\to \M\cong \mathrm{Rees}\M/(t-1)\mathrm{Rees}\M
	\end{equation*}
	sends $m_i$ to $x_i$, so it is immediate that $\epsilon\cdot \mathrm{coker}f=0$ by our choice of $m_i$. It thus remains to show that $f$ is $\epsilon$-strict.
	
	Let $m\in \mathrm{Fil}_jM$ be in the image of $f$, i.e. in the $\A$-submodule generated by the $x_i$.
	
	Viewing $mt^j$ as an element of the $j$th graded piece of $\mathrm{Rees}\M$, it follows that $\epsilon mt ^j$ is in the $\mathrm{Rees}\A$-submodule generated by the $m_i$. But now $\mathrm{Rees}\M$ is a graded module and the $m_i$ are homogeneous, so this implies that $\epsilon mt^j$ is in the image of the natural map
	\begin{equation*}
		\oplus_i \mathrm{Fil}_{j-d_i}\A\to \mathrm{Fil}\M t^j
	\end{equation*} 
	sending $1$ in the $i$th summand to $x_i t^j= t^{j-d_i}\cdot m_i$.
	
	Quotienting out by $(t-1)$ yields $\epsilon m\in \sum \mathrm{Fil}_{j-d_i} \A\cdot x_i=f(\mathrm{Fil}_j \F)$, as required.
\end{proof}

\begin{cor}
	\label{almostres}
	Let $\M$ be an almost finitely generated $\A$-module, endowed with an almost good filtration. Let $\epsilon \in \mathfrak{m}$. Then $\M$ admits an $\epsilon$-strict $\epsilon$-approximation by finite filt-free $\A$-modules, i.e. there exists an $\epsilon$-exact complex
	\begin{equation*}
		\hdots \to \F^{i}\to \F^{i+1}\to \hdots \to \F^0 \to \M \to 0
	\end{equation*}
	such that each $\F^i$ is finite filt-free and each morphism is $\epsilon$-strict.
\end{cor}
\begin{proof}
	Apply the above inductively, noting that at each step, the kernel of the differential $\F^i\to \F^{i+1}$ is an almost finitely generated $\A$-module whose induced subspace filtration is almost good.
\end{proof}

\begin{rmk}
	It follows from Proposition \ref{Reesapproximation} that in this case $\mathrm{Rees}\F^\bullet$ is an $\epsilon^2$-approximation of $\mathrm{Rees}\M$, and $\mathrm{gr}\F^\bullet$ is an $\epsilon^3$-approximation of $\mathrm{gr}\M$.
\end{rmk}

Let $\F$ be a finite filt-free $\A$-module, and let $\N$ be any filtered $\A$-module. Note that we can endow
\begin{equation*}
	\mathrm{Hom}_\A(\F, \N)
\end{equation*}
with a $\mathbb{Z}$-filtration of $R$-modules as follows: a morphism $\phi: \F\to \N$ is in
\begin{equation*}
	\mathrm{Fil}_i\mathrm{Hom}_\A(\F, \N)
\end{equation*}
if $f(\mathrm{Fil}_j\F)\subseteq \mathrm{Fil}_{i+j}\N$ for each $j$. Explicitly, if 
\begin{equation*}
	\F=\oplus_{i=1}^r \A(n_i),
\end{equation*}
then $\mathrm{Hom}_\A(\F, \N)$ can be identified as an $R$-module with $\oplus_{i=1}^r\N$, as any morphism is determined by where the free generators are sent, and as filtered $R$-modules we have
\begin{equation*}
	\mathrm{Hom}_\A(\F, \N)\cong \oplus_{i=1}^r \N(-n_i).
\end{equation*}

\begin{thm}
	\label{grExtvan}
	Let $\M$ be an almost finitely generated $\A$-module endowed with an almost good filtration. Let $j\geq 0$. If we have
	\begin{equation*}
		\mathrm{Ext}^j_{\mathrm{gr}\A}(\mathrm{gr}\M, \mathrm{gr}\A)\cong^a0,
	\end{equation*}
	then
	\begin{equation*}
		\mathrm{Ext}^j_\A(\M, \A)\cong^a 0.
	\end{equation*}
\end{thm}

\begin{proof}
	Let $\M$ be an almost finitely generated $\A$-module endowed with an almost good filtration. Fix $\epsilon\in \mathfrak{m}$ and let $\F^\bullet\to \M$ be an $\epsilon$-strict $\epsilon$-approximation of $\M$ by finite filt-free modules as in Corollary \ref{almostres}.
	
	Consider the complex $\mathrm{Hom}_\A(\F^\bullet, \A)$. By the discussion above, this is a complex of filtered $\A$-modules, and it is straightforward to see (cf. \cite[chapter III, Lemma 2.2.2.(2), Lemma 2.2.3.(2)]{LVO}) that
	\begin{equation*}
		\mathrm{Rees}\mathrm{Hom}_\A(\F^\bullet, \A)\cong \mathrm{Hom}_{\mathrm{Rees}\A} (\mathrm{Rees}\F^\bullet, \mathrm{Rees}\A)
	\end{equation*}
	and
	\begin{equation*}
		\mathrm{gr}\mathrm{Hom}_\A(\F^\bullet, \A)\cong \mathrm{Hom}_{\mathrm{gr}\A} (\mathrm{gr}\F^\bullet, \mathrm{gr}\A).
	\end{equation*}
	
	By Proposition \ref{Reesapproximation}, $\mathrm{Rees}\F^\bullet\to \mathrm{Rees}\M$ is an $\epsilon^2$-approximation. Thus the cokernel of the natural map
	\begin{equation*}
		\mathrm{Ext}^j_{\mathrm{Rees}\A}(\mathrm{Rees}\M, \mathrm{Rees}\A)\to \mathrm{H}^j(\mathrm{Hom}_{\mathrm{Rees}\A}(\mathrm{Rees}\F^\bullet, \mathrm{Rees}\A))
	\end{equation*}
	is annihilated by $\epsilon^{2(j+2)}$ by Lemma \ref{Extapproximation}.
	
	The short exact sequence
	\begin{equation*}
		\begin{xy}
			\xymatrix{0\ar[r]& \mathrm{Rees}\M\ar[r]^t & \mathrm{Rees}\M\ar[r]& \mathrm{gr}\M\ar[r]&0
			}
		\end{xy}
	\end{equation*}
	and the identification
	\begin{equation*}
		\mathrm{Ext}^{j+1}_{\mathrm{Rees}\A}(\mathrm{gr}\M, \mathrm{Rees}\A)\cong \mathrm{Ext}^j_{\mathrm{gr}\A}(\mathrm{gr}\M, \mathrm{gr}\A)\cong^a 0
	\end{equation*}
	as in \cite[Lemma 1.1]{ASZ} shows that multiplication by $t$ is an almost surjective map on $\mathrm{Ext}^j_{\mathrm{Rees}\A}(\mathrm{Rees}\M, \mathrm{Rees}\A)$. It then follows from Lemma \ref{almostsurjcomp} that multiplication by $t^n$ is almost surjective for all $n\geq 1$.
	
	Entirely analogously, $\mathrm{gr}\F^\bullet\to \mathrm{gr}\M$ is an $\epsilon^3$-approximation, so that the natural map
	\begin{equation*}
		\mathrm{Ext}^j_{\mathrm{gr}\A}(\mathrm{gr}\M, \mathrm{gr}\A)\to \mathrm{H}^j(\mathrm{Hom}_{\mathrm{gr}\A}(\mathrm{gr}\F^\bullet, \mathrm{gr}\A))
	\end{equation*}
	has a cokernel annihilated by $\epsilon^{3(j+2)}$. Since $\mathrm{Ext}^j_{\mathrm{gr}\A}(\mathrm{gr}\M, \mathrm{gr}\A)\cong^a 0$, it follows that $\mathrm{H}^j(\mathrm{Hom}_{\mathrm{gr}\A}(\mathrm{gr}\F^\bullet, \mathrm{gr}\A))$ is annihilated by $\epsilon^{3(j+2)+1}$.
	
	The filtration on $\mathrm{Hom}_\A(\F^{-j}, \A)$ induces a subspace filtration on the kernel of the differential, and thus by passing to the quotient we obtain a filtration on $\mathrm{H}^j(\mathrm{Hom}_\A(\F^\bullet, \A))$. 
	
	Note that there is a natural graded surjection
	\begin{equation*}
		\mathrm{H}^j(\mathrm{Hom}_{\mathrm{Rees}\A}(\mathrm{Rees}\F^\bullet, \mathrm{Rees}\A))\to \mathrm{Rees}\mathrm{H}^j(\mathrm{Hom}_\A(\F^\bullet, \A)),
	\end{equation*}
	which on the $i$th graded piece corresponds to the surjection
	\begin{equation*}
		\frac{\mathrm{ker}\mathrm{d}_j\cap \mathrm{Fil}_i}{\mathrm{d}_{j-1}(\mathrm{Fil}_i)}\to \frac{\mathrm{ker}\mathrm{d}_j\cap \mathrm{Fil}_i}{\mathrm{Im}\mathrm{d}_{j-1}\cap \mathrm{Fil}_i},
	\end{equation*} 
	where $\mathrm{d}_j: \mathrm{Hom}_{\A}(\F^{-j}, \A)\to \mathrm{Hom}_{\A}(\F^{-j-1}, \A)$ denotes the $j$th differential, and we have shortened $\mathrm{Fil}_i\mathrm{Hom}_{\A}(\F^{-j}, \A)$ (resp. $\mathrm{Fil}_i\mathrm{Hom}_{\A}(\F^{-j+1}, \A)$) to $\mathrm{Fil}_i$.
	
	We thus have the following diagram of right $\mathrm{Rees}\A$-modules:
	\begin{equation*}
		\begin{xy}
			\xymatrix{
			\mathrm{Ext}^j_{\mathrm{Rees}\A}(\mathrm{Rees}\M, \mathrm{Rees}\A)\ar[d]^{t^n}\ar[r]&\mathrm{H}^j(\mathrm{Hom}_{\mathrm{Rees}\A}(\mathrm{Rees}\F^\bullet, \mathrm{Rees}\A))\ar[d]^{t^n}\ar[r]&\mathrm{Rees}\mathrm{H}^j(\mathrm{Hom}_\A(\F^\bullet, \A))\ar[d]^{t^n}\\
			\mathrm{Ext}^j_{\mathrm{Rees}\A}(\mathrm{Rees}\M, \mathrm{Rees}\A)\ar[r]&\mathrm{H}^j(\mathrm{Hom}_{\mathrm{Rees}\A}(\mathrm{Rees}\F^\bullet, \mathrm{Rees}\A))\ar[r]&\mathrm{Rees}\mathrm{H}^j(\mathrm{Hom}_\A(\F^\bullet, \A)),
			}
		\end{xy}
	\end{equation*}
	
	where the first vertical arrow is almost surjective,  and in each row, the first arrow has its cokernel annihilated by $\epsilon^{2(j+2)}$, and the second arrow is surjective.
	
	Let $\mathrm{Fil}_{\mathrm{min}}$ denote the first non-zero filtered piece of $\mathrm{H}^j(\mathrm{Hom}_\A(\F^\bullet, \A))$. Note that this is a $\mathrm{Fil}_0\A$-submodule of $\mathrm{gr}\mathrm{H}^j(\mathrm{Hom}_\A(\F^\bullet, \A))$. It follows from \cite[chapter III, Corollary 1.1.8]{LVO} that $\mathrm{gr}\mathrm{H}^j(\mathrm{Hom}_\A(\F^\bullet, \A))$ is a subquotient of $\mathrm{H}^j(\mathrm{Hom}_{\mathrm{gr}\A}(\mathrm{gr}\F^\bullet, \mathrm{gr}\A))$, which we have seen above is annihilated by $\epsilon^{3(j+2)+1}$. Therefore $\epsilon^{3(j+2)+1}$ annihilates $\mathrm{Fil}_{\mathrm{min}}$.
	
	Now let $x\in \mathrm{Fil}_{\mathrm{min}+n}\mathrm{H}^j(\mathrm{Hom}_\A(\F^\bullet, \A))$, which we view as an element of the $(\mathrm{min}+n)$th graded piece of the corresponding Rees module. Then by the diagram above, $\epsilon^{2(j+2)}x$ has a preimage in $\mathrm{Ext}^j_{\mathrm{Rees}\A} (\mathrm{Rees}\M, \mathrm{Rees}\A)$, and thus $\epsilon^{2(j+2)+1}x$ is contained in $t^n \mathrm{Rees}\mathrm{H}^j(\mathrm{Hom}_\A(\F^\bullet, \A))$. Therefore, $\epsilon^{2(j+2)+1}x\in \mathrm{Fil}_\mathrm{min}$, and thus $\epsilon^{5(j+2)+2}x=0$.
	
	As $n$ was arbitrary, this shows that $\epsilon^{5(j+2)+2}$ annihilates $\mathrm{H}^j(\mathrm{Hom}_\A(\F^\bullet, \A))$. But now the natural map
	\begin{equation*}
		\mathrm{Ext}^j_\A(\M, \A)\to \mathrm{H}^j(\mathrm{Hom}_\A(\F^\bullet, \A))
	\end{equation*} 
	has kernel annihilated by $\epsilon^{j+1}$ by Lemma \ref{Extapproximation}. Therefore, $\mathrm{Ext}^j_\A(\M, \A)$ is annihilated by $\epsilon^{5(j+2)+2+j+1}=\epsilon^{6j+13}$. But $\epsilon\in \mathfrak{m}$ was arbitrary and $\mathfrak{m}^2=\mathfrak{m}$, so $\mathrm{Ext}^j_\A(\M, \A)$ is indeed annihilated by $\mathfrak{m}$, as required.
\end{proof}

\section{Almost Auslander regular rings}
\subsection{Definition and relation to Auslander regularity on the generic fibre}
Recall the following definition.

\begin{defn}
	Let $A$ be a (not necessarily commutative) Noetherian ring.
	\begin{enumerate}[(i)]
		\item The \textbf{grade} of a (left or right) $A$-module $M$ is
		\begin{equation*}
			j(M)=\mathrm{min}\{i| \mathrm{Ext}^i_A(M, A)\neq 0\},
		\end{equation*}
		and $\infty$ if no such $i$ exists.
		\item We say that a left $A$-module $M$ satisfies the \textbf{Auslander condition} if for any $i\geq 0$ and any (right) $A$-submodule $N\subseteq \mathrm{Ext}^i_A(M, A)$, we have $j(N)\geq i$. 
		
		\noindent We say that $A$ satisfies the Auslander condition (as a ring) if each finitely generated left $A$-module satisfies the Auslander condition. 
		\item We call $A$ \textbf{Auslander--Gorenstein} if it satisfies the Auslander condition and has finite left and right self-injective dimension.
		\item We call $A$ \textbf{Auslander regular} if it is Auslander--Gorenstein and has finite global dimension. 
	\end{enumerate}
\end{defn}

In the light of the previous section, we introduce the following `almost' version.

\begin{defn}
	Let $\A$ be an almost Noetherian $R/\pi$-algebra.
	\begin{enumerate}[(i)]
		\item We say that a left $\A$-module $\M$ satisfies the \textbf{almost Auslander condition} if for any $i\geq 0$ and any submodule $\N\subseteq \mathrm{Ext}^i_\A(\M, \A)$, we have
		\begin{equation*}
			\mathrm{Ext}^j_\A(\N, \A)\cong^a 0
		\end{equation*}
		whenever $j<i$.
		
		\noindent We say that $\A$ satisfies the almost Auslander condition if every finitely generated left $\A$-module satisfies the almost Auslander condition.
		\item We call $\A$ \textbf{almost Auslander--Gorenstein} if it satisfies the almost Auslander condition and there exists some $d\geq 0$ such that
		\begin{equation*}
			\mathrm{Ext}^{j}_\A(\M, \A)\cong^a 0
		\end{equation*}
		for any $j>d$ and any finitely generated $\A$-module $\M$. We call the smallest such $d$ the almost self-injective dimension of $\A$.
		\item We call $\A$ \textbf{almost Auslander regular} if it is almost Auslander--Gorenstein and there exists some $d\geq 0$ such that
		\begin{equation*}
			\mathrm{Ext}^{j}_\A(\M, \N)\cong^a 0
		\end{equation*}
		for any $j>d$, any finitely generated $\A$-module $\M$ and any $\A$-module $\N$ which is flat over $R/\pi$. We call the smallest such $d$ the almost global dimension of $\A$.
	\end{enumerate}
\end{defn}

\begin{rmks}
	\leavevmode
	\begin{enumerate}[(i)]
		\item Alternatively, one might require $\M$ in the definitions above to be only almost finitely generated. It is straightforward to verify that this yields an equivalent notion, by approximating an almost finitely generated $\A$-module by a finitely generated one and passing to the corresponding long exact sequence of Ext groups.
		\item The additional flatness assumption in the definition of almost Auslander regularity is added to exclude torsion phenomena like the well-known
		\begin{equation*}
			\mathrm{Ext}^i_{R/\pi}(R/\epsilon, R/\epsilon)\cong R/\epsilon \ \forall i\geq 0
		\end{equation*}
		for $\epsilon\in \mathfrak{m}$ satisfying $\epsilon^2=\pi$.
		\item As a consequence of the above remark, note that the definition is in some sense `relative', i.e. $\A$ might be almost Auslander regular as an $R/\pi$-algebra, but not as an $R/\pi'$-algebra for some $\pi'\neq \pi\in \mathfrak{m}$.
	\end{enumerate}
\end{rmks}

From now on, let $\A$ denote a $\pi$-adically complete, $\pi$-torsionfree $R$-algebra, and write $A=\A\otimes_RK$. We will now discuss how almost properties of $\A/\pi \A$ can be used to deduce properties of $A$. We begin by relating Ext-groups over $\A$ to those over $\A/\pi\A$.

\begin{lem}
	\label{Extfg}
	Suppose that $\A$ is almost Noetherian. 
	\begin{enumerate}[(i)]
		\item Let $\M$ be an almost finitely generated left $\A$-module. Then for any $j\geq 0$, $\mathrm{Ext}^j_{\A}(\M, \A)$ is an almost finitely generated right $\A$-module.
		\item If $\M$ is a finitely generated, $\pi$-torsionfree left $\A$-module, then $\mathrm{Ext}^j_{\A}(\M, \A)$ is finitely generated.
	\end{enumerate}

\end{lem}
\begin{proof}
	\begin{enumerate}[(i)]
		\item If $\M$ is an almost finitely generated $\A$-module and $\epsilon\in \mathfrak{m}$, pick $\delta\in \mathfrak{m}$ such that $\epsilon\in \delta^{j+2}\mathfrak{m}$. Let $\F^\bullet\to \M$ be a $\delta$-approximation of $\M$ by free finitely generated $\A$-modules as remarked after Lemma \ref{Extapproximation}. 
		
		Then Lemma \ref{Extapproximation} implies that
		\begin{equation*}
			\alpha: \mathrm{Ext}^j_\A(\M, \A)\to \mathrm{H}^j(\mathrm{Hom}_\A(\F^\bullet, \A))
		\end{equation*}
		satisfies $\delta^{j+1}\cdot \mathrm{ker}\alpha=0$. Since $\mathrm{Hom}(\F^{-j}, \A)$ is a finitely generated right module and $\A$ is almost Noetherian, it follows that $\mathrm{Im}\alpha$ is almost finitely generated, so let $m_1, \hdots, m_r\in \mathrm{Ext}^j_\A(\M, \A)$ such that $\delta\cdot \mathrm{Im}\alpha \in \sum_i \alpha(m_i)\cdot \A$.
		
		It follows that for any $m\in \mathrm{Ext}^j_\A(\M, \A)$, there exists some $m'\in \sum_i m_i \cdot \A$ such that $\delta\cdot m-m'\in \mathrm{ker}\alpha$, and thus $\delta^{j+2}m=\delta^{j+1}m'\in \sum_i m_i\A$. Hence we have $\epsilon\cdot \mathrm{Ext}^j_\A(\M, \A)\subseteq \sum m_i\cdot \A$, and $\mathrm{Ext}^j_\A(\M, \A)$ is almost finitely generated.	
		
		\item If $\M$ is even finitely generated and $\pi$-torsionfree, let $f: \A^r\to \M$ be a surjection. It follows from Corollary \ref{ptfisfp} that $\M$ is finitely presented, and $\mathrm{ker}f$ is likewise a finitely generated, $\pi$-torsionfree $\A$-module. Continuing in this way, we obtain a free resolution of $\M$ by free left $\A$-modules of finite rank, and hence $\mathrm{RHom}_{\A}(\M, \A)$ can be represented by a complex $(C^\bullet, \mathrm{d}_\bullet)$ consisting of free right $\A$-modules of finite rank. 
		
		Since the kernel of $\mathrm{d}_j$ is a $\pi$-adically saturated right $\A$-submodule of $\A^r$ for some $r$, it follows again from Corollary \ref{ptfisfp} that $\mathrm{ker}\mathrm{d}_j$ is a finitely generated $\A$-module, and hence $\mathrm{Ext}^j_{\A}(\M, \A)$ is also finitely generated.  
	\end{enumerate}
	
\end{proof}

\begin{lem}
	\label{Extmodp}
	Let $\M$ be an $\A/\pi\A$-module. Then there is a natural isomorphism
	\begin{equation*}
		\mathrm{Ext}^{j+1}_\A(\M, \A)\cong \mathrm{Ext}^j_{\A/\pi\A}(\M, \A/\pi\A)
	\end{equation*}
	for each $j$.
\end{lem}
\begin{proof}
	See \cite[Lemma 1.1]{ASZ}.
\end{proof}

\begin{lem}
	\label{Extvanishmodp}
	Suppose that $\A$ is almost Noetherian. Let $\M$ be an almost finitely generated $\pi$-torsionfree left $\A$-module and $i\geq 0$. If 
	\begin{equation*}
		\mathrm{Ext}^i_{\A/\pi\A}(\M/\pi \M, \A/\pi \A)\cong^a 0,
	\end{equation*}
	then 
	\begin{equation*}
		\mathrm{Ext}^i_\A(\M, \A)\cong^a0.
	\end{equation*}
\end{lem}

\begin{proof}
	The short exact sequence
	\begin{equation*}
		\begin{xy}
			\xymatrix{
			0\ar[r]& \M\ar[r]^{\pi}& \M\ar[r]& \M/\pi \M\ar[r]& 0
			}
		\end{xy}
	\end{equation*}
	yields a long exact sequence
	\begin{equation*}
		\hdots \to \mathrm{Ext}^i_\A(\M/\pi \M, \A)\to \mathrm{Ext}^i_\A(\M, \A)\to \mathrm{Ext}^i_\A(\M, \A)\to \mathrm{Ext}^{i+1}_\A(\M/\pi \M, \A)\to\hdots
	\end{equation*}
	By Lemma \ref{Extmodp}, $\mathrm{Ext}^{i+1}_\A(\M/\pi \M, \A)\cong \mathrm{Ext}^i_{\A/\pi\A}(\M/\pi \M, \A/\pi\A)$, so if this Ext group is almost zero, then multiplication by $\pi$ is an almost surjective map on $\mathrm{Ext}^i_\A(\M, \A)$. Thus 
	\begin{equation*}
		\frac{\mathrm{Ext}^i_\A(\M, \A)}{\pi \mathrm{Ext}^i_\A(\M, \A)}\cong^a 0.
	\end{equation*}
	
	By Lemma \ref{Extfg}.(i), $\mathrm{Ext}^i_\A(\M, \A)$ is almost finitely generated. Hence
	\begin{equation*}
		\mathrm{Ext}^i_\A(\M, \A)\cong^a 0
	\end{equation*} 
	by the same argument as in \cite[Lemma 2.5.19]{BZ}.
\end{proof}

It is straightforward to see that the analogous statement also holds for almost finitely generated $\pi$-torsionfree right $\A$-modules.

\begin{thm}
	\label{liftAuslander}
	Suppose that $\A/\pi\A$ satisfies the almost Auslander condition. Then $A$ satisfies the Auslander condition.
\end{thm}
\begin{proof}
	Note that $\A$ is almost Noetherian and $A$ is a two-sided Noetherian $K$-algebra, by Proposition \ref{liftNoeth}.
	
	Let $M$ be a finitely generated $A$-module, and let $N\subseteq \mathrm{Ext}^i_A(M, A)$ be a right submodule.
	
	Choose a surjection $A^r\to M$ and let $\M$ denote the image of $\A^r$, a finitely generated $\pi$-torsionfree $\A$-module lattice. Note that
	\begin{equation*}
		\mathrm{Ext}^i_A(M, A)\cong \mathrm{Ext}^i_\A(\M, \A)\otimes_R K.
	\end{equation*} 
	Let $E$ denote the image of the natural morphism
	\begin{equation*}
		\theta: \mathrm{Ext}^i_\A(\M, \A)\to \mathrm{Ext}^i_A(M, A),
	\end{equation*}
	so that $E$ is naturally isomorphic to the quotient of $\mathrm{Ext}^i_\A(\M, \A)$ by its torsion submodule.
	
	Let $\N'=N\cap E$, and let $\N\subseteq \mathrm{Ext}^i_\A(\M, \A)$ denote its preimage, i.e.
	\begin{equation*}
		\N=\{m\in \mathrm{Ext}^i_\A(\M, \A): \theta(m)\in N\}.	
	\end{equation*}   
	As $N$ is a $K$-vector space, $\N$ is a saturated right $\A$-submodule of $\mathrm{Ext}^i_\A(\M, \A)$.
	
	As $\M$ is $\pi$-torsionfree, we have a short exact sequence
	\begin{equation*}
		0\to \M\to \M\to \M/\pi \M\to 0,
	\end{equation*}
	where the first map is given by multiplication by $\pi$, and we thus have a long exact sequence of cohomology
	\begin{equation*}
		\hdots \to \mathrm{Ext}^i_\A(\M/\pi \M, \A)\to \mathrm{Ext}^i_\A(\M, \A)\to \mathrm{Ext}^i_\A(\M, \A)\to \mathrm{Ext}^{i+1}_\A(\M/\pi \M, \A)\to\hdots
	\end{equation*}
	
	Since $\N$ is saturated, $\N/\pi \N$ is a submodule of $\frac{\mathrm{Ext}^i_\A(\M, \A)}{\pi\mathrm{Ext}^i_\A(\M, \A)}$, which by the long exact sequence above is in turn a submodule of $\mathrm{Ext}^{i+1}_\A(\M/\pi \M, \A)$.
	
	Hence Lemma \ref{Extmodp} implies that $\N/\pi \N$ is a submodule of $\mathrm{Ext}^i_{\A/\pi\A}(\M/\pi \M, \A/\pi\A)$, so by assumption,
	\begin{equation*}
		\mathrm{Ext}^j_{\A/\pi\A}(\N/\pi \N, \A/\pi\A)\cong^a 0
	\end{equation*}
	for all $j<i$.
	
	Now let $T\subset \N$ be all the torsion elements of $\N$, i.e. 
	\begin{equation*}
		T=\{n\in \N: \pi^s n=0 \ \text{for some } s\in \mathbb{N} \}.
	\end{equation*}
	Since $\N$ is saturated, $T$ actually consists of all the torsion elements in $\mathrm{Ext}^i_\A(\M, \A)$. In particular, $T$ is also a saturated submodule of $\mathrm{Ext}^i_\A(\M, \A)$. Applying the argument above to $T$ instead of $\N$ therefore yields
	\begin{equation*}
		\mathrm{Ext}^j_{\A/\pi\A}(T/\pi T, \A/\pi\A)\cong^a 0
	\end{equation*}
	for all $j<i$.
	
	The short exact sequence
	\begin{equation*}
		0\to T\to \N\to \N'\to 0
	\end{equation*}
	remains exact after tensoring with $R/\pi$, since $\N'$ is $\pi$-torsionfree, and we thus obtain a long exact sequence
	\begin{align*}
		\hdots\to \mathrm{Ext}^i_{\A/\pi\A}(\N'/\pi \N', \A/\pi\A)\to \mathrm{Ext}^i_{\A/\pi\A}(\N/\pi\N, \A/\pi\A)\\
		\to \mathrm{Ext}^i_{\A/\pi\A}(T/\pi T, \A/\pi\A)\to \mathrm{Ext}^{i+1}_{\A/\pi\A}(\N'/\pi\N', \A/\pi\A)\to \hdots
	\end{align*}
	Therefore, we have
	\begin{equation*}
		\mathrm{Ext}^j_{\A/\pi\A}(\N'/\pi\N', \A/\pi\A)\cong^a0
	\end{equation*}
	for all $j<i$.
	
	Finally, as $\N'$ is an almost finitely generated $\pi$-torsionfree $\A$-module, the right module version of Lemma \ref{Extvanishmodp} implies that
	\begin{equation*}
		\mathrm{Ext}^j_\A(\N', \A)\cong^a0
	\end{equation*}
	for all $j<i$, and thus
	\begin{equation*}
		\mathrm{Ext}^j_A(N, A)\cong \mathrm{Ext}^j_\A(\N', \A)\otimes_RK\cong 0,
	\end{equation*}
	for all $j<i$, as required.
\end{proof}

\begin{thm}
	Suppose $\A/\pi\A$ is almost Auslander--Gorenstein of almost self-injective dimension $d$. Then $A$ is Auslander--Gorenstein of self-injective dimension $\leq d$.
\end{thm}

\begin{proof}
	We have already seen in Theorem \ref{liftAuslander} that $A$ satisfies the Auslander condition.
	
	Now suppose that for any finitely generated $\A/\pi\A$-module $\M$ we have
	\begin{equation*}
		\mathrm{Ext}^j_{\A/\pi\A}(\M, \A/\pi\A)\cong^a0
	\end{equation*}
	for $j>d$.
	
	Let $M$ be a finitely generated $A$-module. Let $A^r\to M$ be a surjective $A$-module map, and let $\M$ denote the image of $\A^r$.
	
	As $\M$ is a finitely generated $\A$-module, $\M/\pi\M$ is finitely generated over $\A/\pi\A$. Therefore
	\begin{equation*}
		\mathrm{Ext}^j_{\A/\pi\A}(\M/\pi\M, \A/\pi\A)\cong^a0
	\end{equation*}
	for $j>d$ by assumption. 
	But now Lemma \ref{Extvanishmodp} forces
	\begin{equation*}
		\mathrm{Ext}^j_\A(\M, \A)\cong^a0,
	\end{equation*}
	since $\M$ is finitely generated and $\pi$-torsionfree, and hence
	\begin{equation*}
		\mathrm{Ext}^j_A(M, A)\cong \mathrm{Ext}^j_\A(\M, \A)\otimes_R K=0
	\end{equation*}
	for $j>d$, completing the proof.
\end{proof}

Note that in applying Lemma \ref{Extvanishmodp}, we used the almost finite generation of $\mathrm{Ext}^j_{\A}(\M, \A)$ as a right $\A$-module to ensure that we could apply the almost version of Nakayama's lemma (\cite[Lemma 2.5.19]{BZ}). In order to employ similar techniques for $R$-modules like $\mathrm{Ext}^j_{\A}(\M, \N)$, we need to make a brief detour via derived complete modules.

\begin{defn}
	Let $\epsilon\in \mathfrak{m}$. For an $R$-module $\M$, let $T(M, \epsilon)\in \mathrm{D}(R-\mathrm{mod})$ denote a homotopy limit of the system
	\begin{equation*}
		\begin{xy}
			\xymatrix{\hdots\ar[r]^{\epsilon}&\M\ar[r]^{\epsilon}&\M\ar[r]^{\epsilon}&\M,}
		\end{xy}
	\end{equation*}
	i.e. there exists a distinguished triangle
	\begin{equation*}
		T(\M, \epsilon)\to \prod_n \M\to \prod_n \M,
	\end{equation*}
	where the map $\alpha_{\epsilon}:\prod_n \M\to \prod_n \M$ sends $(x_0, x_1, \hdots)$ to $(x_0-\epsilon x_1, x_1-\epsilon x_2, \hdots)$.
	We say that $\M$ is \textbf{derived $\epsilon$-complete} if $T(\M, \epsilon)\cong 0$, i.e. $\M$ is derived $\epsilon$-complete if and only if $\alpha_{\epsilon}$ is an isomorphism.
	
	A complex of $R$-modules is derived $\epsilon$-complete if each of its cohomology groups is (this is equivalent to the definition given in \cite[tag 091S]{stacksproj} thanks to \cite[tag 091U]{stacksproj}).
\end{defn}

\begin{lem}
	\label{derivedcompl}
	\leavevmode
	\begin{enumerate}[(i)]
		\item An $R$-module is $\epsilon$-adically complete if and only if it is $\epsilon$-adically separated and derived $\epsilon$-adically complete.
		\item Let $\M$ and $\N$ be left $\A$-modules, and assume that $\N$ is derived $\epsilon$-complete. Then $\mathrm{Ext}^j_{\A}(\M, \N)$ is a derived $\epsilon$-complete $R$-module for any $j$.
		\item (Almost version of derived Nakayama.) Suppose that $\M$ is a derived $\epsilon$-adically complete $R$-module for any $\epsilon \in \mathfrak{m}$ and $\M/\pi \M\cong^a 0$. Then $\M\cong^a 0$. 
		\item Let $\M$, $\N$ be left $\A$-modules. If $\N$ is finitely generated and
		\begin{equation*}
			\frac{\mathrm{Ext}^j_{\A}(\M, \N)}{\pi \mathrm{Ext}^j_{\A}(\M, \N)}\cong^a 0,
		\end{equation*} 
		then $\mathrm{Ext}^j_{\A}(\M, \N)\cong^a 0$.
	\end{enumerate}
\end{lem}

\begin{proof}
	\begin{enumerate}[(i)]
		\item See \cite[tag 091T]{stacksproj}.
		\item By \cite[tag 0A6E]{stacksproj}, $\mathrm{RHom}_{\A}(\M, \N)$ is derived $\epsilon$-complete as a complex, so $\mathrm{Ext}^j_{\A}(\M, \N)$ is derived $\epsilon$-adically complete.
		\item This is the almost version of \cite[tag 09B9]{stacksproj}. Fix an arbitrary $\epsilon\in \mathfrak{m}$ with $|\epsilon|>|\pi|$, and let $x_0\in \M$. We wish to show that $\epsilon\cdot x_0=0$. By assumption, $\epsilon\M=\pi \M$, so there exists $x_1\in \M$ such that $\epsilon x_0=\pi x_1=(\pi \epsilon^{-1})\epsilon x_1$. Repeating in this manner, we obtain an element $(\epsilon x_0, \epsilon x_1, \hdots)\in \prod_n \M$ which is in the kernel of the morphism
		\begin{equation*}
			\alpha_{\pi\epsilon^{-1}}: \prod_n \M\to \prod_n\M.
		\end{equation*} 
		Since $\M$ is derived $\pi \epsilon^{-1}$-adically complete, this element is forced to be zero, so in particular $\epsilon\cdot x_0=0$, as required.
		\item By Lemma \ref{almostAR}, $\N$ is $\pi$-adically complete and hence $\epsilon$-adically complete for any $\epsilon\in \mathfrak{m}$. Hence $\mathrm{Ext}^j_{\A}(\M, \N)$ is derived $\epsilon$-adically complete for any $\epsilon\in \mathfrak{m}$ by (ii). We can thus apply (iii) to deduce that $\mathrm{Ext}^j_{\A}(\M, \N)\cong^a 0$. \qedhere
	\end{enumerate}
\end{proof}

\begin{thm}
	\label{liftAreg}
	Suppose that $\A/\pi\A$ is almost Auslander regular of almost global dimension $d$. Then $A$ is Auslander regular of global dimension $\leq d$.
\end{thm}
\begin{proof}
	By the theorem above, it suffices to prove that $\mathrm{gl.dim.}A\leq d$.
	
	Let $M$ and $N$ be finitely generated $A$-modules. As before, pick finitely generated ($\pi$-torsionfree) $\A$-submodules $\M$, $\N$ such that $\M\otimes_RK\cong M$, $\N\otimes_RK\cong N$.
	
	Since $\N$ is $\pi$-torsionfree, $\N/\pi\N$ is flat over $R/\pi$, so that by assumption
	\begin{equation*}
		\mathrm{Ext}^j_{A/\pi\A}(\M/\pi\M, \N/\pi\N)\cong^a0
	\end{equation*}
	for $j>d$.
	
	Since $\M$ is flat over $R$, we have
	\begin{equation*}
		\A/\pi\A\otimes^\mathbb{L}_\A\M\cong R/\pi \otimes^\mathbb{L}_R\M\cong \M/\pi\M,
	\end{equation*}
	and hence
	\begin{equation*}
		\mathrm{R}\mathrm{Hom}_{\A}(\M, \N/\pi\N)\cong \mathrm{R}\mathrm{Hom}_{\A/\pi\A}(\M/\pi\M, \N/\pi\N).
	\end{equation*}
	Therefore
	\begin{equation*}
		\mathrm{Ext}^j_\A(\M, \N/\pi\N)\cong^a0
	\end{equation*}
	for $j>d$.
	
	Now consider the long exact sequence
	\begin{equation*}
		\hdots\to \mathrm{Ext}^i_\A(\M, \N)\to \mathrm{Ext}^i_\A(\M, \N)\to \mathrm{Ext}^i_\A(\M, \N/\pi\N)\to \hdots
	\end{equation*}
	induced by
	\begin{equation*}
		\begin{xy}
			\xymatrix{0\ar[r]&\N\ar[r]^{\pi}& \N\ar[r]&\N/\pi\N\ar[r]&0},
		\end{xy}
	\end{equation*} 
	using the fact that $\N$ is $\pi$-torsionfree.
	
	This shows that for any $j>d$, multiplication by $\pi$ is an almost surjective map on $\mathrm{Ext}^j_\A(\M, \N)$, i.e.
	\begin{equation*}
		\frac{\mathrm{Ext}^j_\A(\M, \N)}{\pi \mathrm{Ext}^j_\A(\M, \N)}\cong^a0.
	\end{equation*}
	
	Since $\N$ is a finitely generated $\A$-module, Lemma \ref{derivedcompl}.(iv) now implies that $\mathrm{Ext}^j_\A(\M, \N)\cong^a 0$, and thus
	\begin{equation*}
		\mathrm{Ext}^j_A(M, N)\cong \mathrm{Ext}^j_\A(\M, \N)\otimes_RK=0
	\end{equation*}
	for $j>d$. Hence $A$ is Auslander regular of global dimension $\leq d$.
\end{proof}
\subsection{Preliminary results on the almost Auslander condition}
\begin{lem}
	Let $\A$ be an $R$-algebra and let
	\begin{equation*}
		0\to \M'\to \M\to \M''\to 0
	\end{equation*}
	be a short exact sequence of $\A$-modules. Let $\epsilon, \delta\in \mathfrak{m}$ and let $j\geq 0$.
	\begin{enumerate}[(i)]
		\item If $\epsilon\cdot \M'=0$ and $\delta\cdot\mathrm{Ext}^j_{\A}(\M, \A)=0$, then $\epsilon\delta\cdot \mathrm{Ext}^j_{\A}(\M'', \A)=0$.
		\item If $\epsilon\cdot \M''=0$ and $\delta\cdot\mathrm{Ext}^j_{\A}(\M, \A)= 0$, then $\epsilon\delta\cdot\mathrm{Ext}^j_{\A}(\M', \A)=0$.
		\item If $\epsilon\cdot \mathrm{Ext}^j_{\A}(\M', \A)=0$ and $\delta\cdot \mathrm{Ext}^j_{\A}(\M'', \A)=0$, then $\epsilon\delta \mathrm{Ext}^j_{\A}(\M, \A)=0$.
	\end{enumerate}
\end{lem}
\begin{proof}
	All of these are immediate from the corresponding long exact sequence of Ext groups.
\end{proof}

\begin{lem}
	\label{sesExta0}
	Let $\N$ be a module over some $R$-algebra $\A$, and let $j\geq 0$. Suppose that for every $\epsilon\in \mathfrak{m}$, there exists an $\epsilon$-exact sequence
	\begin{equation*}
		\begin{xy}
			\xymatrix{
			0\ar[r]&\N'\ar[r]^{\alpha}&\N\ar[r]^{\beta}&\N''\to 0
			}
		\end{xy}
	\end{equation*}
	with $\mathrm{Ext}^j_{\A}(\N', \A)\cong^a \mathrm{Ext}^j_{\A}(\N'', \A)\cong^a 0$.
	
	Then $\mathrm{Ext}^j_{\A}(\N, \A)\cong^a 0$.
\end{lem}
\begin{proof}
	Fix $\epsilon\in \mathfrak{m}$ and consider an $\epsilon$-exact sequence as in the statement of the lemma. We can decompose this into the following short exact sequences:
	\begin{equation*}
		0\to \mathrm{ker}\alpha\to \N'\to \mathrm{Im}\alpha\to 0,
	\end{equation*}
	\begin{equation*}
		0\to \mathrm{Im}\alpha\to \mathrm{ker}\beta\to Q_1\to 0,
	\end{equation*}
	\begin{equation*}
		0\to \mathrm{ker}\beta\to\N\to \mathrm{Im}\beta\to 0,
	\end{equation*}
	and
	\begin{equation*}
		0\to \mathrm{Im}\beta\to \N''\to Q_2\to 0,
	\end{equation*}
	where $\mathrm{ker}\alpha$, $Q_1$ and $Q_2$ are annihilated by $\epsilon$.
	
	Applying repeatedly the previous lemma shows that $\mathrm{Ext}^j_{\A}(\N, \A)$ is annihilated by $\epsilon^3\mathfrak{m}^2$. Explicitly: $\mathrm{Ext}^j_{\A}(\mathrm{Im}\alpha, \A)$ is annihilated by $\epsilon\mathfrak{m}$ by (i), and hence $\mathrm{Ext}^j_{\A}(\mathrm{ker}\beta, \A)$ is annihilated by $\epsilon^2\mathfrak{m}$ by (iii). Since $\mathrm{Ext}^j_{\A}(\mathrm{Im}\beta, \A)$ is annihilated by $\epsilon\mathfrak{m}$ by (ii), it follows that $\mathrm{Ext}^j_{\A}(\N, \A)$ is annihilated by $\epsilon^3\mathfrak{m}^2$ by (iii).

	Since $\epsilon$ was arbitrary and $\mathfrak{m}^2=\mathfrak{m}$, it follows that $\mathrm{Ext}^j_{\A}(\N, \A)$ is almost zero. 
\end{proof}

Let $\epsilon\in \mathfrak{m}$. In the following, we say that an $\A$-module $\M$ satisfies the \textbf{$\epsilon$-Auslander condition in degree $i$} if for any submodule $\N\subseteq \mathrm{Ext}^i_{\A}(\M, \A)$, we have
\begin{equation*}
	\epsilon\cdot \mathrm{Ext}^j_{\A}(\N, \A)=0
\end{equation*}
for any $j<i$.

Note that in this case, we can deduce
\begin{equation*}
	\epsilon^2\mathrm{Ext}^j_{\A}(\N, \A)=0
\end{equation*}
for any $j<i$ and any subquotient $\N$ of $\mathrm{Ext}^i_{\A}(\M, \A)$, by considering the long exact sequence of Ext groups.

\begin{lem}
	\label{epsAus}
	Let $\A$ be an $R$-algebra and let
	\begin{equation*}
		0\to \M'\to \M\to \M''\to 0
	\end{equation*}
	be a short exact sequence of $\A$-modules. Let $\epsilon, \delta\in \mathfrak{m}$, and let $i\geq 0$.
	\begin{enumerate}[(i)]
		\item If $\epsilon\cdot \M'=0$ and $\M$ satisfies the $\delta$-Auslander condition in degree $i$, then $\M''$ satisfies the $\epsilon\cdot\delta$-Auslander condition in degree $i$.
		\item If $\epsilon\cdot \M''=0$ and $\M$ satisfies the $\delta$-Auslander condition in degree $i$, then $\M'$ satisfies the $\epsilon\cdot \delta^2$-Auslander condition in degree $i$.
		\item If $\M'$ satisfies the $\epsilon$-Auslander condition in degree $i$ and $\M''$ satisfies the $\delta$-Auslander condition in degree $i$, then $\M$ satisfies the $\epsilon\cdot \delta^2$-Auslander condition in degree $i$.
	\end{enumerate}
\end{lem}
\begin{proof}
	We show (i), the other parts follow the same strategy. Let $\N\subseteq \mathrm{Ext}^i_{\A}(\M'', \A)$. The long exact sequence of Ext groups
	\begin{equation*}
		\begin{xy}
			\xymatrix{
			\hdots \ar[r]&\mathrm{Ext}^{i-1}_{\A}(\M', \A)\ar[r]^f&\mathrm{Ext}^i_{\A}(\M'', \A)\ar[r]^g& \mathrm{Ext}^i_{\A}(\M, \A)\ar[r]&\hdots
			}
		\end{xy}
	\end{equation*}
	yields a short exact sequence
	\begin{equation*}
		0\to \N\cap \mathrm{ker}g\to \N\to g(\N)\to 0.
	\end{equation*}
	Now $\N\cap \mathrm{ker}g$ is a subquotient of $\mathrm{Ext}^{i-1}_{\A}(\M', \A)$, and hence is annihilated by $\epsilon$, and $g(\N)$ is a submodule of $\mathrm{Ext}^i_{\A}(\M, \A)$. Thus the result follows from the corresponding long exact sequence of Ext groups.
\end{proof}

\begin{lem}
	\label{sesforaAcond}
	Let $\M$ be a module over some $R$-algebra $\A$. Suppose that for every $\epsilon\in \mathfrak{m}$, there exists an $\epsilon$-exact sequence
	\begin{equation*}
		\begin{xy}
			\xymatrix{
			0\ar[r]&\M'\ar[r]^{\alpha}&\M\ar[r]^{\beta}&\M''\ar[r]& 0
			}
		\end{xy}
	\end{equation*}
	such that $\M'$ and $\M''$ satisfy the almost Auslander condition. Then $\M$ satisfies the almost Auslander condition.
\end{lem}
\begin{proof}
	We wish to show that $\M$ satisfies the $\epsilon$-Auslander condition in degree $i$ for each $\epsilon\in \mathfrak{m}$ and each $i\geq 0$.
	
	Fix $\epsilon\in \mathfrak{m}$ and consider an $\epsilon$-exact sequence as in the statement of the lemma. Decomposing into short exact sequences as in Lemma \ref{sesExta0}, it follows from repeatedly applying Lemma \ref{epsAus} that $\M$ satisfies the $\epsilon^6$-Auslander condition in degree $i$, for any $i\geq 0$.
	
	Since $\epsilon$ was arbitrary and $\mathfrak{m}^2=\mathfrak{m}$, the result follows.
\end{proof}

\section{Almost Auslander regularity of polynomial rings}

Let $\A$ be an almost Auslander regular $R/\pi$-algebra, and suppose that the polynomial ring $\A[X, Y]$ is almost Noetherian.

We endow $\A[x]$ with the usual filtration given by the degree of $x$, so that $\mathrm{gr}\A[x]\cong \A[x]$, $\mathrm{Rees}\A[x]\cong \A[X, Y]$ (with $X$ corresponding to $xt$ and $Y$ corresponding to $t$). In particular, this filtration is almost Zariskian by assumption. The goal of this section is to adapt the arguments in \cite[chapter III, sections 2.2, 2.3]{LVO} to the almost setting, showing that $\A[x]$ is almost Auslander regular.

\subsection{Reduction to graded modules}

We first have the following reduction step.
\begin{thm}
	\label{reducetogr}
	Let $\M$ be an almost finitely generated $\A[x]$-module, equipped with an almost good filtration. If $\mathrm{gr}\M$ satisfies the almost Auslander condition, then so does $\M$.
\end{thm}

\begin{proof}
	Let $\N\subseteq \mathrm{Ext}^i_{\A[x]}(\M, \A[x])$ be a right $\A[x]$-submodule, and let $\epsilon\in \mathfrak{m}$.
	
	Choose an $\epsilon$-strict $\epsilon$-approximation $\F^\bullet\to \M$ as in Corollary \ref{almostres} and consider the natural morphism
	\begin{equation*}
		\alpha: \mathrm{Ext}^i_{\A[x]}(\M, \A[x])\to \mathrm{H}^i(\mathrm{Hom}_{\A[x]}(\F^\bullet, \A[x])).
	\end{equation*}
	Let $\N'=\alpha(\N)$, so that we have a short exact sequence
	\begin{equation*}
		0\to \N\cap \mathrm{ker}\alpha\to \N\to \N'\to 0.
	\end{equation*}
	
	We abbreviate $E^i=\mathrm{Ext}^i_{\mathrm{Rees}\A[x]}(\mathrm{Rees}\M, \mathrm{Rees}\A[x])$. As in the proof of Theorem \ref{grExtvan}, we have an exact sequence
	\begin{equation*}
		\begin{xy}
			\xymatrix{E^i\ar[r]^t& E^i\ar[r]&\mathrm{Ext}^i_{\mathrm{gr}\A[x]}(\mathrm{gr}\M, \mathrm{gr}\A[x])}
		\end{xy}
	\end{equation*}
	Since $\mathrm{gr}\M$ satisfies the almost Auslander condition, it follows that whenever $T\subseteq E^i$, then
	\begin{equation*}
		\mathrm{Ext}^j_{\mathrm{gr}\A[x]}(\frac{T}{T\cap tE^i}, \mathrm{gr}\A[x])\cong^a 0
	\end{equation*}
	for any $j<i$.
	
	Recall from the proof of Theorem \ref{grExtvan} that the cokernel of the natural map
	\begin{equation*}
		 E^i\to \mathrm{Rees}\mathrm{H}^i(\mathrm{Hom}_{\A[x]}(\F^\bullet, \A[x]))
	\end{equation*}
	is annihilated by $\epsilon^{2(i+2)}$, and thus the same is true of the natural map
	\begin{equation*}
		\beta: E^i\to \mathrm{gr}\mathrm{H}^i(\mathrm{Hom}_{\A[x]}(\F^\bullet, \A[x])),
	\end{equation*}
	which is obtained from the above by passing to the quotient.
	 
	If we equip $\N'\subseteq \mathrm{H}^i(\mathrm{Hom}_{\A[x]}(\F^\bullet, \A[x]))$ with the subspace filtration, we have
	\begin{equation*}
		\mathrm{gr}\N'\subseteq \mathrm{gr}\mathrm{H}^i(\mathrm{Hom}_{\A[x]}(\F^\bullet, \A[x])),
	\end{equation*}
	and it follows that $\epsilon^{2i+4}\mathrm{gr}\N'$ is a subquotient of $E^i$.
	
	Let $T'\subset T\subset E^i$ such that 
	\begin{equation*}
		T/T'\cong \epsilon^{2i+4}\mathrm{gr}\N'.
	\end{equation*}
	Since $T\cap tE^i$ is contained in the kernel of $\beta$, this yields a short exact sequence
	\begin{equation*}
		0\to \frac{T'}{T'\cap tE^i}\to \frac{T}{T\cap tE^i}\to \epsilon^{2i+4}\mathrm{gr}\N'\to 0,
	\end{equation*}
	and we can deduce from the above that
	\begin{equation*}
		\mathrm{Ext}^j_{\mathrm{gr}\A[x]}(\epsilon^{2i+4}\mathrm{gr}\N', \mathrm{gr}\A[x])\cong^a 0
	\end{equation*}
	for all $j<i$.
	
	Therefore, multiplication by $t$ is almost surjective on the module $\mathrm{Ext}^j_{\mathrm{Rees}\A[x]}(\epsilon^{2i+4}\mathrm{Rees}\N', \mathrm{Rees}\A[x])$ for all $j<i$, by considering the corresponding long exact sequence of Ext groups for
	\begin{equation*}
		\begin{xy}
			\xymatrix{0\ar[r]&\epsilon^{2i+4}\mathrm{Rees}\N'\ar[r]^t&\epsilon^{2i+4}\mathrm{Rees}\N'\ar[r]&\epsilon^{2i+4}\mathrm{gr}\N'\ar[r]&0}
		\end{xy}
	\end{equation*}
	and noting again that
	\begin{equation*}
		\mathrm{Ext}^{j+1}_{\mathrm{Rees}\A[x]}(\epsilon^{2i+4}\mathrm{gr}\N', \mathrm{Rees}\A[x])\cong \mathrm{Ext}^j_{\mathrm{gr}\A[x]}(\epsilon^{2i+4}\mathrm{gr}\N', \mathrm{gr}\A[x]).
	\end{equation*}
	This implies that multiplication by $t^n$ is an almost surjective map on $\mathrm{Ext}^j_{\mathrm{Rees}A[x]}(\epsilon^{2i+4}\mathrm{Rees}\N', \mathrm{Rees}\A[x])$ for any $j<i$ and any $n\geq 0$.
	
	Since multiplication by $\epsilon^{2i+4}$ on $\mathrm{Ext}^j_{\mathrm{Rees}\A[x]}(\mathrm{Rees}\N', \mathrm{Rees}\A[x])$ factors through $\mathrm{Ext}^j_{\mathrm{Rees}\A[x]}(\epsilon^{2i+4}\mathrm{Rees}\N', \mathrm{Rees}\A[x])$, it follows that 
	\begin{equation*}
		\mathfrak{m}\cdot \epsilon^{2i+4}\mathrm{Ext}^j_{\mathrm{Rees}\A[x]}(\mathrm{Rees}\N', \mathrm{Rees}\A[x])\subseteq t^n\mathrm{Ext}^j_{\mathrm{Rees}\A[x]}(\mathrm{Rees}\N', \mathrm{Rees}\A[x])
	\end{equation*}
	for any $n\geq 0$.
	
	Since the filtration on $\N'$ is again an almost good filtration, we can choose an $\epsilon$-strict $\epsilon$-approximation $\G^\bullet\to \N'$ by finite filt-free $\A[x]$-modules, yielding an $\epsilon^2$-approximation $\mathrm{Rees}\G^\bullet\to \mathrm{Rees}\N'$ by finite free $\mathrm{Rees}\A[x]$-modules by Proposition \ref{Reesapproximation}.
	
	As in Theorem \ref{grExtvan}, we have the commutative diagram
	\begin{equation*}
		\begin{xy}
			\xymatrix{\mathrm{Ext}^j_{\mathrm{Rees}\A[x]}(\mathrm{Rees}\N', \mathrm{Rees}\A[x])\ar[d]\ar[r]^{t^n}& \mathrm{Ext}^j_{\mathrm{Rees}\A[x]}(\mathrm{Rees}\N', \mathrm{Rees}\A[x])\ar[d]\\
				\mathrm{H}^j(\mathrm{Hom}_{\mathrm{Rees}\A[x]}(\mathrm{Rees}\G^\bullet, \mathrm{Rees}\A[x]))\ar[r]\ar[d]&\mathrm{H}^j(\mathrm{Hom}_{\mathrm{Rees}\A[x]}(\mathrm{Rees}\G^\bullet, \mathrm{Rees}\A[x]))\ar[d]\\
				\mathrm{Rees}\mathrm{H}^j(\mathrm{Hom}_{\A[x]}(\G^\bullet, \A[x]))\ar[r]& \mathrm{Rees}\mathrm{H}^j(\mathrm{Hom}_{\A[x]}(\G^\bullet, \A[x]))}
		\end{xy}
	\end{equation*}
where in each column, the top arrow has its cokernel annihilated by $\epsilon^{2j+4}$ and the bottom arrow is almost surjective.

Let $\mathrm{Fil}_\mathrm{min}$ denote the first non-trivial filtered piece of $\mathrm{H}^j(\mathrm{Hom}(\G^\bullet, \A[x]))$. It then follows from the above that
\begin{equation*}
	\epsilon^{2i+5}\cdot\epsilon^{2j+4}\mathrm{Rees}\mathrm{H}^j(\mathrm{Hom}(\G^\bullet, \A[x]))\subseteq t \mathrm{Rees}\mathrm{H}^j(\mathrm{Hom}(\G^\bullet, \A[x])),
\end{equation*}
so that $\mathrm{gr}\mathrm{H}^j(\mathrm{Hom}(\G^\bullet, \A[x]))$, and hence $\mathrm{Fil}_\mathrm{min}$, is annihilated by $\epsilon^{2i+2j+9}$.

But now if $x\in \mathrm{Fil}_{\mathrm{min}+n}\mathrm{H}^j(\mathrm{Hom}(\G^\bullet, \A[x]))$, viewed as the $(\mathrm{min}+n)$th graded piece of $\mathrm{Rees}\mathrm{H}^j(\mathrm{Hom}(\G^\bullet, \A[x]))$, then $\epsilon^{2i+2j+9}x$ is contained in $t^n\mathrm{Rees}\mathrm{H}^j(\mathrm{Hom}(\G^\bullet, \A[x]))$ and hence in $t^n\mathrm{Fil}_\mathrm{min}$.

Therefore $\mathrm{Rees}\mathrm{H}^j(\mathrm{Hom}(\G^\bullet, \A[x]))$, and thus $\mathrm{H}^j(\mathrm{Hom}(\G^\bullet, \A[x]))$, is annihilated by $\epsilon^{4i+4j+18}$.

Since the natural map
\begin{equation*}
	\mathrm{Ext}^j_{\A[x]}(\N', \A[x])\to \mathrm{H}^j(\mathrm{Hom}_{\A[x]}(\G^\bullet, \A[x]))
\end{equation*}
has kernel annihilated by $\epsilon^{j+1}$ by Lemma \ref{Extapproximation}, it follows that $\mathrm{Ext}^j_{\A[x]}(\N', \A[x])$ is annihilated by $\epsilon^{4i+5j+19}$.

Finally, the long exact sequence of Ext groups corresponding to
\begin{equation*}
	0\to \N\cap \mathrm{ker}\alpha\to \N\to \N'\to 0
\end{equation*}
and the fact that $\mathrm{ker}\alpha$ is annihilated by $\epsilon^{i+1}$ shows that $\mathrm{Ext}^j_{\A[x]}(\N, \A[x])$ is annihilated by $\epsilon^{5(i+j)+20}$. Since $\epsilon$ was arbitrary, this proves that $\mathrm{Ext}^j_{\A[x]}(\N, \A[x])$ is almost zero for any $j< i$, and thus $\M$ satisfies the almost Auslander condition.
\end{proof}

\subsection{The almost Auslander condition for graded modules}
Now let $\M$ be a finitely generated $\A[x]$-module. The quotient filtration induced by a surjection
\begin{equation*}
	\A[x]^r\to \M
\end{equation*}
endows $\M$ with an almost good filtration (as $\mathrm{Rees}\M$ is a quotient of $\mathrm{Rees}\A[x]^r$).

It follows that $\mathrm{gr}\M$ is an almost finitely generated graded $\A[x]$-module such that each graded piece $\mathrm{gr}_i\M$ is a finitely generated $\A$-module. In fact, $\mathrm{gr}\M$ is even finitely generated over $\A[x]$.

In the light of Theorem \ref{reducetogr}, we now wish to show that $\mathrm{gr}\M$ satisfies the almost Auslander condition.
\begin{lem}
	\label{sesforgr}
	For any $\epsilon\in \mathfrak{m}$, there exists an exact sequence of graded $\A[x]$-modules
	\begin{equation*}
		\A[x]\otimes_\A\M'\to \mathrm{gr}\M\to \M''\to 0,
	\end{equation*}
	where
	\begin{enumerate}[(i)]
		\item $\M'$ is a finitely generated $\A$-module.
		\item $\M''$ is finitely generated as an $\A$-module, and there exists some $d\in \mathbb{N}$ such that $x^d\cdot \M''=0$.
		\item the kernel of the first arrow is annihilated by $\epsilon$.
	\end{enumerate}
\end{lem} 
\begin{proof}
	Fix $\epsilon\in \mathfrak{m}$.
	
	Let 
	\begin{equation*}
		T=\{m\in \mathrm{gr}\M: x^s\cdot m=0\ \text{for some} \ s\in \mathbb{N}\},
	\end{equation*}
	a graded submodule of $\mathrm{gr}\M$. In particular, $T$ is almost finitely generated over $\A[x]$, so let $T^\epsilon$ denote a finitely generated $\A[x]$-submodule of $T$ such that $T/T^\epsilon$ is annihilated by $\epsilon$. Since $T$ is graded, we can assume without loss of generality that $T^\epsilon$ is generated by homogeneous elements $m_1, \hdots, m_r\in T$, of degree $d_1, \hdots, d_r$, say. By definition, we can find $s_i\in \mathbb{N}$ such that $x^{s_i}m_i=0$ for each $i$. In particular, $\mathrm{gr}_d T^\epsilon=0$ whenever $d\geq\mathrm{max}\{d_i+s_i\}$, since $\mathrm{gr}_dT^\epsilon\subseteq \sum_i x^{d-d_i}\A m_i$. This in turn implies that $\epsilon\cdot \mathrm{gr}_dT=0$ whenever $d\geq \mathrm{max}\{d_i+s_i\}$.
	
	Since $\mathrm{gr}\M$ is finitely generated and graded, it is generated by some homogeneous elements $n_1, \hdots, n_{r'}$, of degrees $e_1, \hdots, e_{r'}$, say. In particular, if $d\geq \mathrm{max}\{e_i\}$, then the natural map 
	\begin{equation*}
		\A[x]\otimes_\A \mathrm{gr}_d\M\to \oplus_{j\geq d}\mathrm{gr}_j\M
	\end{equation*}
	is surjective: for any $n\in \mathrm{gr}_j\M$ with $j\geq d$, we have 
	\begin{equation*}
		n=\sum_i x^{j-d_i}a_in_i=x^{j-d}\cdot (\sum_i x^{d-d_i}a_in_i)
	\end{equation*}
	for some $a_i\in \A$.
	
	Now pick $d\geq \mathrm{max}\{d_i+s_i, e_j\}$. The natural graded map
	\begin{equation*}
		\theta: \A[x]\otimes_\A \mathrm{gr}_d\M\to \mathrm{gr}\M
	\end{equation*}
	has image $\oplus_{j\geq d}\M$ by the above, so $\mathrm{coker}\theta\cong \oplus_{j<d}\mathrm{gr}_j\M$ as an $\A$-module and $x^d\mathrm{coker}\theta=0$.
	
	Moreover, $\mathrm{ker}\theta$ is contained in $\oplus_{j\geq d} x^{j-d}\otimes\mathrm{gr}_d T$, so $\epsilon\cdot \mathrm{ker}\theta=0$ by the above. Thus we obtain the result with $\M'=\mathrm{gr}_d\M$ and $\M''=\mathrm{coker}\theta=\mathrm{gr}\M/\oplus_{j\geq d}\mathrm{gr}_j\M$.
\end{proof}

We will also need the following variant.
\begin{lem}
	\label{sesExta0gr}
	Let $\N$ be an almost finitely generated graded $\A[x]$-module which is $x$-torsionfree. For any $\epsilon\in \mathfrak{m}$, there exists a sequence of graded $\A[x]$-modules
	\begin{equation*}
		\A[x]\otimes_\A\N'\to \N\to \N'',
	\end{equation*}
	where
	\begin{enumerate}[(i)]
		\item $\N'$ is an almost finitely generated $\A$-module.
		\item $\N''$ is an almost finitely generated $\A[x]$-module which is already almost finitely generated over $\A$, and there exists some $d\in \mathbb{N}$ such that $x^d\cdot \N''=0$.
		\item $\alpha$ is injective, $\beta$ is surjective, and $\epsilon\cdot \mathrm{ker}\beta\subseteq \mathrm{Im}\alpha\subseteq \mathrm{ker}\beta$. In particular, the above yields a short $\epsilon$-exact sequence.
	\end{enumerate}
\end{lem}
\begin{proof}
	As $\N$ is almost finitely generated and graded, there exists a finitely generated graded submodule $T$ such that $\epsilon\cdot \N\subset T$. We can assume without loss of generality that $T$ is generated by homogeneous elements $n_1, \hdots, n_r\in T$ of degree $d_1, \hdots, d_r$, say.
	
	Let $d\geq \mathrm{max}\{d_i\}$ and consider the natural graded morphism
	\begin{equation*}
		\alpha: \A[x]\otimes_\A\mathrm{gr}_dT\to \N.
	\end{equation*}
	Since $\N$ is $x$-torsionfree, $\alpha$ is injective.
	
	Now let
	\begin{equation*}
		\N''=\N/\oplus_{j\geq d}\mathrm{gr}_j\N
	\end{equation*}
	and let $\beta:\N\to \N''$ denote the natural quotient map. Then $\N''$ is an almost finitely generated $\A[x]$-module, and $x^d\cdot \N''=0$. In particular, $\N''$ is also almost finitely generated over $\A$.
	
	It is clear from construction that $\mathrm{Im}\alpha\subseteq \mathrm{ker}\beta=\oplus_{i\geq d} \mathrm{gr}_i\N$.
	
	It remains to show that $\epsilon\cdot \mathrm{ker}\beta\subseteq \mathrm{Im}\alpha$. If $n\in \mathrm{gr}_i\N$ with $i\geq d$, then $\epsilon n\in \mathrm{gr}_i T$, and thus
	\begin{equation*}
		\epsilon\cdot n=\sum_i x^{i-d_i}a_in_i=x^{i-d}\cdot (\sum_i x^{d-d_i}a_in_i)
	\end{equation*}  
	for some $a_i\in \A$. Thus $\epsilon\mathrm{gr}_i\N\subseteq \mathrm{Im}\alpha$, as required.
\end{proof}

\begin{lem}
	\label{slicegr}
	Let $\N$ be a graded $\A[x]$-module, $\N=\oplus_{i\geq 0}\mathrm{gr}_i\N$, and suppose that 
	\begin{equation*}
		\mathrm{Ext}^j_\A(\mathrm{gr}_s\N, \A)\cong^a 0
	\end{equation*}
	for each $s$.
	Then
	\begin{equation*}
		\mathrm{Ext}^{j+1}_{\A[x]}(\frac{\N}{\oplus_{j\geq d}\mathrm{gr}_j\N}, \A[x])\cong^a 0
	\end{equation*}
	for each $d$.
\end{lem}
\begin{proof}
	Let $T_s$ denote the $\A[x]$-module which is the $\A$-module $\mathrm{gr}_s\N$ with $x$ acting trivially.
	Then
	\begin{equation*}
		\mathrm{Ext}^{j+1}_{\A[x]}(T_s, \A[x])\cong \mathrm{Ext}^j_\A(\mathrm{gr}_s\N, \A)\cong^a 0
	\end{equation*}
	by assumption.
	
	Now argue inductively, using the short exact sequences of $\A[x]$-modules
	\begin{equation*}
		0\to T_{d-1}\to \frac{\N}{\oplus_{j\geq d}\N}\to \frac{\N}{\oplus_{j\geq d-1}\N}\to 0. \qedhere
	\end{equation*}
\end{proof}

Looking at Lemma \ref{sesforgr} and Lemma \ref{sesforaAcond}, it is now clear that we only have to consider two special cases to establish the Auslander condition for finitely generated graded modules. The first case is quite straightforward:
\begin{lem}
	\label{finitecase}
	Let $\M''$ be an $\A[x]$-module which is already almost finitely generated over $\A$, satisfying $x^d\cdot \M''=0$ for some $d\in \mathbb{N}$. Then $\M''$ satisfies the almost Auslander condition.
\end{lem}
\begin{proof}
	Suppose first that we can take $d=1$, i.e. $x\cdot \M''=0$. Then 
	\begin{equation*}
		\mathrm{Ext}^i_{\A[x]}(\M'', \A[x])\cong \mathrm{Ext}^{i-1}_\A(\M'', \A)
	\end{equation*}
	by \cite[Lemma 1.1]{ASZ}, and any $\N\subseteq \mathrm{Ext}^i_{\A[x]}(\M'', \A[x])$ is also an $\A[x]$-module satisfying $x\cdot \N=0$.
	Thus
	\begin{equation*}
		\mathrm{Ext}^j_{\A[x]}(\N, \A[x])\cong \mathrm{Ext}^{j-1}_\A(\N, \A)\cong^a0
	\end{equation*} 
	for any $j<i$, as $\A$ is almost Auslander regular.
	
	For $d>1$, consider the short exact sequence
	\begin{equation*}
		0\to x^{d-1}\cdot \M''\to \M''\to \M''/x^{d-1}\M''\to 0.
	\end{equation*}
	Each term is an $\A[x]$-module which is almost finitely generated over $\A$, since $\A$ is almost Noetherian. We can thus use induction and Lemma \ref{sesforaAcond} to finish the proof.
\end{proof}

It remains to consider the case of a graded module of the form $\A[x]\otimes \M'$ for some finitely generated $\A$-module $\M'$.

First note that for each $i$,
\begin{equation*}
	E^i:=\mathrm{Ext}^i_{\A[x]}(\A[x]\otimes_\A\M', \A[x])\cong \mathrm{Ext}^i_\A(\M', \A)\otimes_\A \A[x]
\end{equation*}
carries a natural grading. If $\N\subseteq E^i$ is a right submodule, the grading induces a natural filtration on $\N$, and $\mathrm{gr}\N$ is still an $\A[x]$-submodule of $E^i$, with $\mathrm{gr}_n\N$ being an $\A$-submodule of
\begin{equation*}
	\mathrm{Ext}^i_\A(\M', \A)\otimes x^n\cong \mathrm{Ext}^i_\A(\M', \A).
\end{equation*}

The filtration on $E^i\cong \mathrm{Ext}^i_{\A}(\M', \A)\otimes_{\A}\A[x]$ is easy to see to be an almost good filtration, since $\mathrm{Ext}^i_{\A}(\M', \A)$ is almost finitely generated over $\A$ by the same argument as in Lemma \ref{Extfg}.(i). It follows that the subspace filtration on $\N$ is also an almost good filtration.

\begin{prop}
	\label{gradedsub}
	Let $\M'$ be a finitely generated $\A$-module. If $\N$ is a graded submodule of $\mathrm{Ext}^i_{\A[x]}(\A[x]\otimes_\A\M', \A[x])$, then
	\begin{equation*}
		\mathrm{Ext}^j_{\A[x]}(\N, \A[x])\cong^a0
	\end{equation*}
	for any $j<i$.
\end{prop}
\begin{proof}
	We wish to apply Lemma \ref{sesExta0} and Lemma \ref{sesExta0gr}.
	
	By the same argument as in Lemma \ref{Extfg}.(i), $E^i$ is almost finitely generated, so $\N$ is an almost finitely generated graded $\A[x]$-module, which is $x$-torsionfree, as $E^i$ is. Fix $\epsilon\in \mathfrak{m}$. By Lemma \ref{sesExta0gr}, there now exists an $\A$-submodule $\N'\subseteq \mathrm{gr}_d\N$ for some $d$ such that 
	\begin{equation*}
		\alpha: \A[x]\otimes_\A \N'\to \N
	\end{equation*}
	is an injection, and $\epsilon\cdot (\oplus_{j\geq d}\mathrm{gr}_j\N)\subseteq \mathrm{Im}\alpha$.
	
	As $\N'\subseteq \mathrm{gr}_d\N\subseteq \mathrm{Ext}^i_\A(\M', \A)\otimes x^d\cong \mathrm{Ext}^i_{\A}(\M', \A)$, the Auslander regularity of $\A$ yields
	\begin{equation*}
		\mathrm{Ext}^j_{\A[x]}(\A[x]\otimes_\A\N', \A[x])\cong \mathrm{Ext}^j_\A(\N', \A)\otimes_\A \A[x]\cong^a 0
	\end{equation*}
	for all $j<i$.
	
	Moreover, for any $s$, $\mathrm{gr}_s\N$ is an $\A$-submodule of $\mathrm{Ext}^i_\A(\M', \A)\otimes x^s$ for any $s$, so
	\begin{equation*}
		\mathrm{Ext}^j_\A(\mathrm{gr}_s\N, \A)\cong^a 0
	\end{equation*} 
	for any $j<i$.
	
	Therefore, Lemma \ref{slicegr} implies that
	\begin{equation*}
		\mathrm{Ext}^j_{\A[x]}(\frac{\N}{\oplus_{j\geq d}\mathrm{gr}_j\N}, \A[x])\cong^a 0.
	\end{equation*}
	
	Now we can apply Lemma \ref{sesExta0} to obtain that 
	\begin{equation*}
		\mathrm{Ext}^j_{\A[x]}(\N, \A[x])\cong^a0
	\end{equation*}
	for any $j< i$, as required.
\end{proof}

\begin{thm}
	Let $\A$ be an almost Auslander regular $R/\pi$-algebra, and assume that $\A[X, Y]$ is almost Noetherian. Then $\A[x]$ satisfies the almost Auslander condition.
\end{thm}
\begin{proof}
	Let $\M$ be a finitely generated $\A[x]$-module. Endow it with an almost good filtration via a surjection $\A[x]^r\to \M$. By Theorem \ref{reducetogr}, it suffices to show that $\mathrm{gr}\M$ satisfies the almost Auslander condition.
	
	By Lemma \ref{sesforgr} and Lemma \ref{sesforaAcond}, we only need to consider two cases:
	
	If $\mathrm{gr}\M$ is finitely generated over $\A$, with $x^d\mathrm{gr}\M=0$ for some $d$, then $\mathrm{gr}\M$ satisfies the almost Auslander condition by Lemma \ref{finitecase}.
	
	If $\mathrm{gr}\M=\A[x]\otimes_\A \M'$ for some finitely generated $\A$-module and $\N\subseteq \mathrm{Ext}^i_{\A[x]}(\mathrm{gr}\M, \A[x])$, then we have seen above that $\N$ is endowed with an almost good filtration such that $\mathrm{gr}\N\subseteq \mathrm{Ext}^i_{\A[x]}(\mathrm{gr}\M, \A[x])$. Proposition \ref{gradedsub} now implies that
	\begin{equation*}
		\mathrm{Ext}^j_{\A[x]}(\mathrm{gr}\N, \A[x])\cong^a0
	\end{equation*}
	for any $j<i$. Now we can invoke Theorem \ref{grExtvan} to deduce that
	\begin{equation*}
		\mathrm{Ext}^j_{\A[x]}(\N, \A[x])\cong^a0,
	\end{equation*}
	showing that $\mathrm{gr}\M$ satisfies the almost Auslander condition.
	
	Thus every finitely generated $\A[x]$-module satisfies the almost Auslander condition, so $\A[x]$ satisfies the almost Auslander condition.
\end{proof}
\subsection{Dimension results}

As before, let $\A$ be an almost Auslander regular $R/\pi$-algebra such that $\A[X, Y]$ is almost Noetherian. Suppose that $\A$ has almost global dimension $d$ in the sense that
\begin{equation*}
	\mathrm{Ext}^j_\A(\M, \N)\cong^a0
\end{equation*}
for any $j>d$, whenever $\M$ is a finitely generated $\A$-module and $\N$ is an $\A$-module which is flat over $R/\pi$.

Now let $\M$ be a finitely generated $\A[x]$-module. As before, we equip $\M$ with an almost good filtration via some surjection $\A[x]^r\to \M$.

We will relate $\M$ to $\A[x]\otimes_{\A}\mathrm{Fil}_s\M$ by studying the map
\begin{align*}
	\theta: \M\to\M\\
	m\mapsto xm.
\end{align*}

\begin{lem}
	For any $\epsilon\in \mathfrak{m}$, there exists some integer $s$ such that 
	\begin{equation*}
		\begin{xy}
			\xymatrix{(\mathrm{Fil}_{s-1}\M\ar[r]^{\theta}& \mathrm{Fil}_s\M)\ar[r]&(\M\ar[r]^{\theta}&\M)}
		\end{xy}\tag{$\star$}
	\end{equation*}
	is an $\epsilon$-approximation, i.e. $\epsilon\cdot \mathrm{ker}\theta\subseteq \mathrm{ker}\theta\cap \mathrm{Fil}_{s-1}\M$ and $\epsilon\M\subseteq \mathrm{Fil}_s\M+x\M$.
\end{lem}
\begin{proof}
	Let $\delta\in \mathfrak{m}$ such that $\epsilon\in \delta^2\mathfrak{m}$.
	Since $\M$ is almost finitely generated, there exists some $s$ such that the cokernel of the natural morphism
	\begin{equation*}
		\A[x]\otimes_\A\mathrm{Fil}_{s-1}\M\to \M
	\end{equation*}
	is annihilated by $\delta$, and $\delta\cdot \{m\in \M| xm=0\}\subseteq \mathrm{Fil}_{s-1}\M$.
	
	It follows immediately that the natural morphisms
	\begin{equation*}
		\mathrm{ker}(\mathrm{Fil}_{s-1}\M\to \mathrm{Fil}_s\M)\to \mathrm{ker}(\M\to \M)
	\end{equation*}
	and 
	\begin{equation*}
		\mathrm{coker}(\mathrm{Fil}_{s-1}\M\to \mathrm{Fil}_s\M)\to \mathrm{coker}(\M\to \M)
	\end{equation*}
	have kernel and cokernel annihilated by $\delta$, so the cone of $(\star)$ is $\epsilon$-exact.
\end{proof}

Now note that for any $n$, we have the following commutative diagram
\begin{equation*}
	\begin{xy}
		\xymatrix{0\ar[r]& \A[x]\otimes_\A \mathrm{Fil}_{n-1}\M\ar[r]\ar[d]^{\alpha_n}& \A[x]\otimes_\A \M\ar[r]\ar[d]^{\beta}& \A[x]\otimes_\A \frac{\M}{\mathrm{Fil}_{n-1}\M}\ar[d]^{\gamma_n}\ar[r]&0\\
		0\ar[r]& \A[x]\otimes_\A \mathrm{Fil}_{n}\M\ar[r]& \A[x]\otimes_\A \M\ar[r]& \A[x]\otimes_\A \frac{\M}{\mathrm{Fil}_n\M}\ar[r]&0
		}
	\end{xy}
\end{equation*}
where each row is exact and the vertical maps are all induced by $x\otimes \mathrm{id}-\mathrm{id}\otimes \theta$.

Note that $\beta$ is injective and $\mathrm{coker}\beta\cong \M$ (see e.g. the proof of \cite[Theorem 4.3.7]{Weibel}). We thus have an exact sequence
\begin{equation*}
	0\to \mathrm{ker}\gamma_n\to \mathrm{coker}\alpha_n\to \M\to \mathrm{coker}\gamma_n\to 0.
\end{equation*}
In particular, since $\M$ and $\mathrm{coker}\alpha_n$ are almost finitely generated, it follows that $\mathrm{ker}\gamma_n$ and $\mathrm{coker}\gamma_n$ are almost finitely generated $\A[x]$-modules.

\begin{lem}
	For any $\epsilon\in \mathfrak{m}$, there exists an integer $s$ such that kernel and cokernel of $\gamma_s$ are annihilated by $\epsilon$.
\end{lem}
\begin{proof}
	For any $n$, consider the commutative diagram
	\begin{equation*}
		\begin{xy}
			\xymatrix{0\ar[r]& \A[x]\otimes_\A\mathrm{gr}_{n}\M\ar[r]\ar[d]& \A[x]\otimes_\A\frac{\M}{\mathrm{Fil}_{n-1}\M}\ar[r]\ar[d]& \A[x]\otimes_\A \frac{\M}{\mathrm{Fil}_n\M}\ar[r]\ar[d]& 0\\
				0\ar[r]& \A[x]\otimes_\A\mathrm{gr}_{n+1}\M\ar[r]& \A[x]\otimes_\A\frac{\M}{\mathrm{Fil}_n\M}\ar[r]& \A[x]\otimes_\A \frac{\M}{\mathrm{Fil}_{n+1}\M}\ar[r]& 0}
		\end{xy}
	\end{equation*}
	with exact rows and each vertical arrow induced by $x\otimes \mathrm{id}-\mathrm{id}\otimes x$. In particular, the first vertical arrow is simply multiplication by $-x$ on the second factor, which is surjective by construction of the filtration on $\M$.
	
	It follows that both $\mathrm{ker}\gamma_n\to \mathrm{ker}\gamma_{n+1}$ and $\mathrm{coker}\gamma_n\to \mathrm{coker}\gamma_{n+1}$ are surjective.
	
	Since $\mathrm{ker}\gamma_1$ and $\mathrm{coker}\gamma_1$ are almost finitely generated, there exists some $s$ such that
	\begin{equation*}
		\epsilon \cdot \mathrm{ker}\gamma_1\subseteq \A[x]\otimes_\A \frac{\mathrm{Fil}_s\M}{\mathrm{Fil}_1\M}
	\end{equation*}
	and 
	\begin{equation*}
		\epsilon\cdot \mathrm{coker}\gamma_1\subseteq \frac{(\A[x]\otimes_\A \mathrm{Fil}_s\M)+\mathrm{Im}\gamma_1}{\mathrm{Im}\gamma_1}.
	\end{equation*}
	But then $\epsilon\cdot \mathrm{ker}\gamma_s\subseteq \A[x]\otimes_\A\frac{\mathrm{Fil}_s\M}{\mathrm{Fil}_s\M}=0$ and $\epsilon\cdot \mathrm{coker}\gamma_s=0$, as required.
\end{proof}

\begin{prop}
	For any $\epsilon\in \mathfrak{m}$, there exists some integer $s$ and an $\epsilon$-exact sequence
	\begin{equation*}
		0\to \A[x]\otimes_\A \mathrm{Fil}_{s-1}\M\to \A[x]\otimes_\A \mathrm{Fil}_s\M\to \M\to 0. 
	\end{equation*}
\end{prop}
\begin{proof}
	Recall that we have a short exact sequence
	\begin{equation*}
		0\to \mathrm{ker}\gamma_n\to \mathrm{coker}\alpha_n\to \M\to \mathrm{coker}\gamma_n\to 0
	\end{equation*}
	for any integer $n$. Choosing $s$ as in the previous lemma, this shows that
	\begin{equation*}
		\begin{xy}
			\xymatrix{0\ar[r]&\A[x]\otimes_\A\mathrm{Fil}_{s-1}\M\ar[r]^{\alpha_s}&\A[x]\otimes_\A \mathrm{Fil}_s\M\ar[r]&\M\ar[r]&0}
		\end{xy}
	\end{equation*}
	is $\epsilon$-exact. Explicitly, we know that $\alpha_s$ is injective, $\epsilon\cdot \M$ is contained in the image of $\mathrm{coker}\alpha_s$, since $\epsilon\cdot \mathrm{coker}\gamma_s=0$, and if $m\in \A[x]\otimes_\A \mathrm{Fil}_s\M$ gets mapped to zero in $\M$, then $\overline{m}\in \mathrm{coker}\alpha_s$ is annihilated by $\epsilon$, as it is in the image of $\mathrm{ker}\gamma_s\to \mathrm{coker}\alpha_s$, and thus $\epsilon m\in \mathrm{Im}\alpha_s$.
\end{proof}
\begin{lem}
	\label{sesfordim}
	Let $\epsilon\in \mathfrak{m}$, and let
	\begin{equation*}
		\begin{xy}
			\xymatrix{0\ar[r]&\M_1\ar[r]^{\alpha}&\M_2\ar[r]^{\beta}&\M_3\ar[r]&0}
		\end{xy}
	\end{equation*}
	be an $\epsilon$-exact sequence of $\A[x]$-modules. Let $\N$ be some left $\A[x]$-module. If there exists $d\geq 0$ such that 
	\begin{equation*}
		\mathrm{Ext}^j_{\A[x]}(\M_1, \N)\cong^a \mathrm{Ext}^j_{\A[x]}(\M_2, \N)\cong^a 0
	\end{equation*}
	for all $j>d$, then 
	\begin{equation*}
		\epsilon^3\mathfrak{m}\mathrm{Ext}^j_{\A[x]}(\M_3, \N)=0
	\end{equation*}
	for all $j>d+1$.
\end{lem}
\begin{proof}
	Decompose the above into short exact sequences
	\begin{equation*}
		0\to \mathrm{ker}\alpha\to \M_1\to \mathrm{Im}\alpha\to 0,
	\end{equation*}
	\begin{equation*}
		0\to \mathrm{Im}\alpha\to \mathrm{ker}\beta\to Q_1\to 0,
	\end{equation*}
	\begin{equation*}
		0\to \mathrm{ker}\beta\to \M_2\to \mathrm{Im}\beta\to 0
	\end{equation*}
	and
	\begin{equation*}
		0\to \mathrm{Im}\beta\to \M_3\to Q_2\to 0,
	\end{equation*}
	where $\mathrm{ker}\alpha, Q_1$ and $Q_2$ are annihilated by $\epsilon$. Working along the corresponding long exact sequences of Ext groups straightforwardly yields the result.
\end{proof}
\begin{thm}
	\label{aAregdimThm}
	Let $\A$ be an almost Auslander regular $\R/\pi$-algebra of almost dimension $d$, and suppose that $\A[x, y]$ is almost Noetherian. Then $\A[x]$ is almost Auslander regular of almost dimension $d+1$.
\end{thm}
\begin{proof}
	We have already seen that $\A[x]$ satisfies the almost Auslander condition, so it suffices to prove the statement on the almost global dimension.
	
	Let $\epsilon\in \mathfrak{m}$.
	
	Let $\M$ be a finitely generated $\A[x]$-module, with an almost good filtration as above, and consider the $\epsilon$-exact sequence
	
	\begin{equation*}
		0\to \A[x]\otimes_\A\mathrm{Fil}_{s-1}\M\to \A[x]\otimes_\A\mathrm{Fil}_s\M\to \M\to 0
	\end{equation*}
	from above. 
	Note that by construction, $\mathrm{Fil}_{s-1}\M$ and $\mathrm{Fil}_s\M$ are finitely generated $\A$-modules.
	
	Let $\N$ be an $\A[x]$-module which is flat over $R/\pi$. Since
	\begin{equation*}
		\mathrm{Ext}^j_{\A[x]}(\A[x]\otimes_A \mathrm{Fil}_s\M, \N)\cong \mathrm{Ext}^j_\A(\mathrm{Fil}_s\M, \N)\cong^a 0
	\end{equation*}
	for any $j>d$, likewise for $\mathrm{Fil}_{s-1}\M$, Lemma \ref{sesfordim} yields that
	\begin{equation*}
		\epsilon^4\mathrm{Ext}^j_{\A[x]}(\M, \N)=0.
	\end{equation*}
	for any $j>d+1$. Since $\epsilon$ was arbitrary, this proves that
	\begin{equation*}
		\mathrm{Ext}^j_{\A[x]}(\M, \N)\cong^a 0
	\end{equation*}
	for any $j>d+1$.
	
	Moreover, if there exists a finitely generated $\A$-module $\M'$ and an $\A$-module $\N'$ which is flat over $R/\pi$ such that
	\begin{equation*}
		\mathrm{Ext}^d_\A(\M', \N')
	\end{equation*}
	is not almost zero, then consider $\M'$ as an $\A[x]$-module with the trivial $x$-action to obtain that
	\begin{equation*}
		\mathrm{Ext}^{d+1}_{\A[x]}(\M', \A[x]\otimes_A\N')\cong \mathrm{Ext}^d_\A(\M', \frac{\A[x]\otimes\N'}{x\A[x]\otimes\N'})\cong \mathrm{Ext}^d_\A(\M', \N') 
	\end{equation*}
	is not almost zero.
\end{proof}

As a consequence, we can now verify that the polynomial ring $R/\pi[x_1, \hdots, x_m]$ is almost Auslander regular, proving Theorem \ref{polyaAreg} from the introduction. 

We argue by induction on $m$. The base case $m=0$ is the following Lemma.
\begin{lem}
	\leavevmode
	\begin{enumerate}[(i)]
		\item Let $\M$ be a finitely generated $R/\pi$-module. We have $\mathrm{Ext}_{R/\pi}^i(\M, R/\pi)\cong^a 0$ for all $i>0$. In particular, $R/\pi$ is almost Auslander--Gorenstein.
		\item Suppose that $\M$ is a finitely generated $R/\pi$-module, and $\N$ is a flat $R/\pi$-module. Then $\mathrm{Ext}^i_{R/\pi}(\M, \N)\cong^a0$ for all $i>0$. In particular, $R/\pi$ is almost Auslander regular, of almost global dimension $0$. 
	\end{enumerate}	
\end{lem}

\begin{proof}
	\leavevmode
	\begin{enumerate}[(i)]
		\item Note that $\M$ can be regarded as a finitely generated $R$-module, and
		\begin{equation*}
			\mathrm{Ext}_{R/\pi}^i(\M, R/\pi)\cong \mathrm{Ext}^{i+1}_R(\M, R)
		\end{equation*}
		by \cite[Lemma 1.1]{ASZ}. 
		For any $\epsilon\in \mathfrak{m}$, the almost Noetherianity of $R$ implies that there exists a short exact sequence
		\begin{equation*}
			0\to S\to \M'\to \M\to 0,
		\end{equation*}
		where $\M'$ is a finitely presented $R$-module and $S$ is annihilated by $\epsilon$.
		By \cite[tag 0ASP]{stacksproj}, $\M'$ is a direct summand of a direct sum of modules of the form $R/f$ for various $f\in R$. In particular, $\mathrm{Ext}^i_R(\M', R)=0$ for all $i>1$. The long exact sequence of Ext groups then implies that $\mathrm{Ext}^i_R(\M, R)$ is annihilated by $\epsilon$ whenever $i>1$. Since $\epsilon$ was arbitrary, this forces $\mathrm{Ext}^i_R(\M, R)\cong^a0$. Thus $\mathrm{Ext}^i_{R/\pi}(\M, R/\pi)\cong^a0$ for all $i>0$, and $\M$ satisfies the almost Auslander condition for $R/\pi$-modules trivially.
		\item This follows directly from \cite[Proposition 2.6.20]{BZ} and the above, noting that any almost Noetherian ring is almost coherent (\cite[Lemma 2.7.7]{BZ}) and any flat module is almost flat (\cite[Lemma 2.2.6]{BZ}). \qedhere
	\end{enumerate}
\end{proof}

\begin{cor}
	The polynomial algebra $R/\pi[x_1, \hdots, x_m]$ is almost Auslander regular of almost global dimension $m$.
\end{cor}
\begin{proof}
	This is now immediate: The lemma above treats the case $m=0$, and Theorem \ref{aAregdimThm} allows us to finish the proof by induction.
\end{proof}

\begin{cor}
	\label{polyreduction}
	Let $\A$ be a $\pi$-adically complete, $\pi$-torsionfree $R$-algebra with $\A/\pi \A\cong R/\pi [x_1, \hdots, x_d]$ for some $d\geq 0$. Then $\A\otimes_R K$ is Auslander regular of global dimension $\leq d$.
\end{cor}
\begin{proof}
	Combine the corollary above with Theorem \ref{liftAreg}.
\end{proof}

\section{Examples: Banach algebras with polynomial reduction}

We can now return to the examples from the introduction, proving Theorem \ref{MainThm1}. 

\subsection{Completed enveloping algebras}
If $\mathfrak{g}$ is a finite-dimensional Lie algebra over $K$ and $\mathfrak{g}_R$ is an $R$-Lie lattice, then $[\pi^n\mathfrak{g}_R, \pi^n\mathfrak{g}_R]\subseteq \pi\cdot \pi^n\mathfrak{g}_R$ for any $n\geq 1$, so that the PBW theorem implies that the $R$-algebra
\begin{equation*}
	U_n^\circ=\h{U_R(\pi^n\mathfrak{g}_R)}
\end{equation*}
satisfies
\begin{equation*}
	U_n^\circ/\pi U_n^\circ \cong \frac{U_R(\pi^n\mathfrak{g}_R)}{\pi U_R(\pi^n\mathfrak{g}_R)}\cong R/\pi [x_1, \hdots, x_d],
\end{equation*}
for $n\geq 1$, where $d=\mathrm{dim}_K \mathfrak{g}$. As an $R$-module, $U_R(\pi^n\mathfrak{g}_R)$ is isomorphic to a polynomial ring over $R$ in $d$ variables, so that its $\pi$-adic completion $U_n^\circ$ is also $\pi$-torsionfree. We can thus apply Corollary \ref{polyreduction} to deduce that $U_n=U_n^\circ\otimes_RK$ is Auslander regular of global dimension $\leq d$.

If we consider the trivial $U_n$-module $K$, computing $\mathrm{R}\mathrm{Hom}_{U_n}(K, U_n)$ via a suitable Koszul complex shows that the global dimension is precisely $d$.

In fact, there is a free resolution of $K$ as a left $U(\mathfrak{g})$-module given by
\begin{equation*}
	0\to U(\mathfrak{g})\otimes_K\wedge^d\mathfrak{g}\to \hdots\to U(\mathfrak{g})\otimes_K \wedge^0\mathfrak{g}\to K,
\end{equation*}
where the map $U(\mathfrak{g})\otimes\wedge^i\mathfrak{g}\to U(\mathfrak{g})\otimes_K\wedge^{i-1}\mathfrak{g}$ is 
\begin{align*}
	f\otimes a_1\wedge a_2\wedge\hdots \wedge a_i\mapsto \sum_j (-1)^jfa_j\otimes a_1\wedge \hdots \wedge \h{a_j}\wedge\hdots\wedge a_i\\
	+\sum_{j, j'} (-1)^{j+j'}f\otimes [a_j, a_{j'}]\wedge a_1\wedge \hdots \wedge \h{a_j} \wedge \hdots \wedge \h{a_{j'}}\wedge\hdots \wedge a_i.
\end{align*}
We can thus compute that $\mathrm{RHom}_{U(\mathfrak{g})}(K, U(\mathfrak{g}))\cong K[-d]$.

It follows from \cite[Theorem 2.13]{SixOp} that $U(\mathfrak{g})\to \w{U(\mathfrak{g})}$ is flat, and since $\w{U(\mathfrak{g})}$ is a Fr\'echet--Stein algebra by \cite[Theorem 2.8]{SixOp}, the map $\w{U(\mathfrak{g})}\to U_n$ is flat by \cite[Remark 3.2]{ST}. Since $U_n\otimes_{U(\mathfrak{g})}K\cong K$, this forces
\begin{equation*}
	\mathrm{Ext}^j_{U_n}(K, U_n)\cong \mathrm{Ext}^j_{U(\mathfrak{g})}(K, U(\mathfrak{g}))\otimes_{U(\mathfrak{g})}U_n
\end{equation*} 
for all $j$, so $\mathrm{RHom}_{U_n}(K, U_n)\cong K[-d]$, as required.
\subsection{Completed Weyl algebras}
Let
\begin{equation*}
	\A=R\langle x_1, \hdots, x_m, \pi^n\frac{\mathrm{d}}{\mathrm{d}x_1}, \hdots, \pi^n\frac{\mathrm{d}}{\mathrm{d}x_m}\rangle
\end{equation*}
denote the unit ball of the completed Weyl algebra $\h{A_{m, n}}$, i.e. $\A$ is the $\pi$-adic completion of the enveloping algebra
\begin{equation*}
	U_{R\langle x_1, \hdots, x_m\rangle}(\pi^n \L)
\end{equation*}
for the $(R\langle x_1, \hdots, x_m\rangle, R)$-Lie--Rinehart algebra $\L=\mathrm{Der}_R(R\langle x_1, \hdots, x_m\rangle)$, cf. \cite{DcapOne}.

We have
\begin{equation*}
	[\pi^n\frac{\mathrm{d}}{\mathrm{d}x_i}, x_j]=\pi^n \delta_{ij}, \ [\pi^n\L, \pi^n\L]\subseteq  \pi^{2n}\L
\end{equation*}
where $\delta_{ij}$ denotes the Kronecker delta, so the commutators vanish modulo $\pi$ as soon as $n\geq 1$.

Thus
\begin{equation*}
	\A/\pi \A\cong U(\pi^n \L)/\pi U(\pi^n\L)\cong R/\pi [X_1, \hdots, X_m, Y_1, \hdots, Y_m]
\end{equation*}
is a polynomial algebra, and $\h{A_{m, n}}$ is Auslander regular of global dimension $\leq 2m$ for any $n\geq 1$ by Corollary \ref{polyreduction}. 

In the appendix, we prove that the global dimension of $\h{A_{m, n}}$ is in fact $m$, as is the case when $K$ is discretely valued (cf. \cite[Corollary 7.4, Theorem 3.3]{AWannals}). Thanks to the Auslander regularity of $\h{A_{m, n}}$, this statement is equivalent to an analogue of Bernstein's inequality over $K$. 

\subsection{Norm-completed distribution algebras}
Lastly, let $G$ denote a compact locally $L$-analytic group of dimension $d$, for some finite extension $L$ of $\mathbb{Q}_p$. Recall from \cite{ST} that the algebra of locally $L$-analytic, $K$-valued distributions $D(G, K)$ admits Banach completions $D_r(G, K)$ for each $1/p\leq r<1$, with $D(G, K)\cong \varprojlim_{r\to 1} D_r(G, K)$.

We first assume that $G$ is an $L$-uniform group as in \cite[Remark 2.2.5]{OS}. Let $r_n=p^{-1/p^n}$. It follows from \cite[Corollary 6.5.4]{Ardakoveq} that
\begin{equation*}
	D_{1/p}(G^{p^n}, K)\cong \h{U_R(p^n\L)}\otimes_R K,
\end{equation*}
where $\L$ is the $R$-Lie algebra associated with $G$ -- in particular, $D_{1/p}(G^{p^n}, K)$ is Auslander regular of dimension $d$ for any $n\geq 1$ by subsection 5.1.

By \cite[Theorem 6.5.11]{Ardakoveq}, $D_{r_n}(G, K)$ contains $D_{1/p}(G^{p^n}, K)$ as a subring such that it is a free left and right $D_{1/p}(G^{p^n}, K)$-module, with a free set of generators given by any set of coset representatives for $G/G^{p^n}$. 

Thus \cite[Corollary 7.3]{Schmidtdistr} implies that $D_{r_n}(G, K)$ is Auslander regular of dimension $d$ for any $n\geq 1$.

Finally, if $G$ is arbitrary, then it admits an open normal $L$-uniform subgroup $H$ by \cite[Lemma 2.2.4]{OS}. Then $D(G, K)$ is a finite free left and right $D(H, K)$-module, with a free set of generators given by any set of coset representatives for $G/H$. For any $1/p\leq r<1$, the Banach algebra $D_r(G, K)$ is then obtained by completing $D(G, K)$ with respect to the semi-norm induced from $D_r(H,K)$ (it follows from the proof of \cite[Theorem 5.1]{ST} that this is indeed an algebra semi-norm). In particular, $D_{r_n}(G, K)$ is a finite free left and right $D_{r_n}(H, K)$-module, with a free set of generators given by coset representatives, and we can apply once more \cite[Corollary 7.3]{Schmidtdistr} to deduce that $D_{r_n}(G, K)$ is Auslander regular of dimension $d=\mathrm{dim}_L G$.

\appendix
\section{Bernstein's inequality for $\h{A_{m, n}}$}	

In this appendix, we adapt the proof of Bernstein's inequality in \cite[Corollary 7.4]{AWannals} to the non-discrete setting, proving the following theorem. 

\begin{thm}
	\label{Bernstein}
	The Tate-Weyl algebra $\h{A_{m, n}}$ satisfies $\mathrm{gl.dim}(\h{A_{m,n}})=m$ for any $n\geq 1$.
	
	In particular, if $M$ is a finitely generated left $\h{A_{m,n}}$-module such that 
	\begin{equation*}
		\mathrm{Ext}^j_{\h{A_{m,n}}}(M, \h{A_{m,n}})=0
	\end{equation*}
	for all $j\leq m$, then $M=0$. In other words: $j(M)\leq m$ for any non-zero finitely generated module $M$.
\end{thm}

In essence, the proof is the same as in \cite{AWannals}, using slightly more algebraic language and a version of Noether normalization over $R/\pi$.

While the proof in \cite{AWannals} expresses Bernstein's inequality in terms of the dimension of the characteristic variety associated to a module, we will remain in our purely algebraic setting. No serious attempt has been made to formally introduce characteristic varieties in the non-discrete setting, even though the arguments below might indicate one possible approach to do so.

\begin{lem}
	\label{injbasechange}
	Let $A$ be an Auslander regular Banach $K$-algebra, $K'$ be a complete field extension of $K$, and set $A'=A\h{\otimes}_K K'$. Then
	\begin{enumerate}[(i)]
		\item $\mathrm{gl.dim.}(A)=\mathrm{inj.dim.}(A)$.
		\item $\mathrm{inj.dim.}(A)\leq \mathrm{inj.dim}(A')$.
		\item $\mathrm{inj.dim.}(A)=\mathrm{inj.dim.}(A')$ if $K'/K$ is finite.
	\end{enumerate}
\end{lem}
\begin{proof}
	As $A$ is Auslander regular, we have $\mathrm{gl.dim.}(A)=\mathrm{inj.dim.}(A)$: if $M, N$ are finitely generated left $A$-modules, then the natural morphism $\mathrm{RHom}_A(M, A)\otimes_A^{\mathbb{L}}N\to \mathrm{RHom}_A(M, N)$ is an isomorphism, since $M$ is perfect. In particular, if $\mathrm{RHom}_A(M, A)$ is concentrated in degrees $\leq d$, then so is $\mathrm{RHom}_A(M, N)$ for all finitely generated $A$-modules $N$. Hence $\mathrm{gl.dim.}(A)\leq \mathrm{inj.dim.}(A)$, and the reverse inequality is obviously true. This proves (i).
	
	If $M$ is a finitely generated $A$-module, then $\mathrm{RHom}_A(M, A)$ can be represented by a bounded strict complex of finite free right $A$-modules and $\mathrm{Ext}^i_A(M, A)$ is a finitely generated right $A$-module, endowed with its canonical Banach structure. Note that
	\begin{equation*}
		\mathrm{Ext}^i_{A'}(K'\h{\otimes}_KM, A')\cong \mathrm{Ext}^i_{A}(M, A')\cong \mathrm{Ext}^i_A(M, A)\h{\otimes}_K K',
	\end{equation*}
	since $-\h{\otimes}_K K'$ is exact on Banach spaces. If $V$ is any Banach $K$-vector space such that $V\h{\otimes}_K K'=0$, then $V=0$ (by \cite[Theorem 10.3.13]{Schikhof}, $V\otimes_K K'$ embeds into $V\h{\otimes}_K K'$), so the above shows that $\mathrm{inj.dim.}(A)\leq \mathrm{inj.dim.}(A')$. Thus (ii) holds. 
	
	For the converse if $K'/K$ is finite, we spell out the argument from \cite[Lemma 4.5]{AWannals}, \cite[Lemma 5.4]{Brown}: let $e_1, \hdots, e_d\in K'$ be a basis of $K'$ over $K$ which is orthogonal with respect to the non-degenerate symmetric bilinear form induced by the trace morphism $\mathrm{tr}: K'\to K$. The trace map now induces a morphism of left $A$-modules $\mathrm{tr}: A'\to A$ such that $\mathrm{tr}(e_iae_j)=a\mathrm{tr}(e_i^2)\delta_{ij}$ for any $a\in A$, since the $e_i$ (or to be precise, the elements $1\otimes e_i$) are central in $A'$. The morphism
	\begin{align*}
		A'\to \mathrm{Hom}_A(A',A)\\
		a\mapsto \left(x\mapsto \mathrm{tr}(xa)\right)
	\end{align*}
	can now easily be checked to be a well-defined morphism of left $A'$-modules, which is in fact an isomorphism with inverse $f\mapsto \sum_{i=1}^d f(e_i)(\mathrm{tr}(e_i^2))^{-1}e_i$. Hence $A'\cong \mathrm{Hom}_A(A', A)\cong \mathrm{R}\mathrm{Hom}_A(A', A)$ as left $A'$-modules and tensor-hom adjunction yields
	\begin{equation*}
		\mathrm{RHom}_A(M',A)\cong \mathrm{RHom}_{A'}(M', \mathrm{RHom}_A(A', A))\cong \mathrm{RHom}_{A'}(M', A')
	\end{equation*}
	for any finitely generated left $A'$-module $M'$. Thus $\mathrm{inj.dim.}(A)\geq \mathrm{inj.dim.}(A')$, proving the result.
\end{proof}

Fix $m\geq 0$ and $n\geq 1$. We will first suppose that the residue field $k$ is algebraically closed. Replacing $K$ by a finite field extension, we can assume that $\pi$ has a square root in $K$, and replacing $\pi$ by its square root, we thus consider the algebra $\h{A_{m, 2n}}$ for $n\geq 1$.

We have performed this modification because we can express $\h{A_{m, 2n}}$ as a complete deformation of another algebra, constructed via a symplectic form.

For this, let $V$ be a free $R$-module of rank $2m$, endowed with the standard non-degenerate anti-symmetric bilinear form $\omega$, i.e. $V$ has a basis $x_1, \hdots, x_m, y_1, \hdots, y_m$ such that
\begin{equation*}
	\omega(x_i, x_j)=\omega(y_i, y_j)=0, \ \omega(y_i, x_j)=\delta_{ij}.
\end{equation*}

As in \cite[Example 3.4.(d)]{AWannals}, we can now define the deformed algebra $K\langle V\rangle_n$ as follows: 

Let $R_{\omega}[V]$ denote the quotient of the tensor $R$-algebra $T_R(V)=\oplus_i V^{\otimes i}$ by the two-sided ideal generated by
\begin{equation*}
	v\otimes w-w\otimes v-\omega(v, w)
\end{equation*}
for $v, w\in V$.

As in \cite[Example 3.5]{AWannals}, the natural filtration on $R_\omega[V]$ makes $R_\omega[V]$ a deformable $R$-algebra in the sense that $\mathrm{gr}R_\omega[V]$ is $\pi$-torsionfree. 

We let $R_\omega[V]_n=\sum \pi^{ni} \mathrm{Fil}_iR_\omega[V]\subseteq R_\omega[V]$ denote the $n$th deformation. This is an $R$-algebra generated by $\pi^nx_1, \hdots, \pi^nx_m, \pi^ny_1, \hdots, \pi^ny_m$ satisfying
\begin{equation*}
	[\pi^nx_i, \pi^nx_j]=[\pi^ny_i, \pi^ny_j]=0, \ [\pi^ny_i, \pi^nx_j]=\pi^{2n}\delta_{ij}.
\end{equation*}
In particular,
\begin{equation*}
	\h{R_\omega[V]_n}\otimes_RK\cong \h{A_{m, 2n}}.
\end{equation*}

In fact, $R_\omega[V]_n$ is naturally isomorphic to $R[x_1,\hdots, x_m, \pi^{2n}\partial_1, \hdots, \pi^{2n}\partial_m]$ by sending $\pi^nx_i$ to $x_i$ and $\pi^ny_i$ to $\pi^{2n}\partial_i$.

Since $n\geq 1$, it follows that 
\begin{equation*}
	\h{R_\omega[V]_n}/\pi\cong R/\pi[X_1, \hdots, X_m, Y_1, \hdots, Y_m]\cong \mathrm{Sym}_{R/\pi}(V/\pi)
\end{equation*}
is a commutative $R/\pi$-algebra.

Let $\mathrm{Gr}_s(V)$ denote the set of all direct summands of $V$ which are free $R$-modules of rank $s$, analogously $\mathrm{Gr}_s(V/\pi)$ and $\mathrm{Gr}_s(V_k)$, where we write $V_k:=V\otimes_Rk$.

We remark that if $W\in \mathrm{Gr}_s(V)$, then the same construction as above yields a ring $\h{R_\omega[W]_n}$, which is naturally a subring of $\h{R_\omega[V]_n}$ such that we have a natural isomorphism
\begin{equation*}
	\h{R_\omega[W]_n}/\pi\cong \mathrm{Sym}_{R/\pi}(W/\pi).
\end{equation*}

In particular, it follows from Corollary \ref{polyreduction} that $\h{R_\omega[W]_n}$ is an almost Auslander regular $R$-algebra.

We use the notation $K\langle W\rangle_n$ to denote the Noetherian Banach algebra $\h{R_\omega[W]_n}\otimes_RK$. In particular, as we have noted above, $\h{A_{m, 2n}}\cong K\langle V\rangle_n$.

We now collect several facts on the Grassmannians $\mathrm{Gr}_s(V)$. 

\begin{lem}
	The natural maps $\mathrm{Gr}_s(V)\to \mathrm{Gr}_s(V/\pi)$ and $\mathrm{Gr}_s(V/\pi)\to \mathrm{Gr}_s(V_k)$ are surjective.
\end{lem} 
\begin{proof}
	Let $W\in \mathrm{Gr}_s(V_k)$ and let $A\in \mathrm{Mat}_d(k)$ be an idempotent matrix whose image is $W$. 
	
	Let $\widetilde{A}\in \mathrm{Mat}_d(R)$ be a $d\times d$-Matrix whose reduction mod $\mathfrak{m}$ is $A$. In particular, $\widetilde{A}^2-I_d\in \mathfrak{m}\mathrm{Mat}_d(R)$, so there exists some $\epsilon\in \mathfrak{m}$ such that $\widetilde{A}$ is idempotent modulo $\epsilon$ (choose $\epsilon\in \mathfrak{m}$ bigger than each matrix entry of $\widetilde{A}^2-I_d$). 
	
	But then it follows from the idempotent lifting lemma \cite[Theorem 21.31]{Lam} that there exists an idempotent matrix $B\in \mathrm{Mat}_d(R)$ such that $\widetilde{A}-B\in \epsilon \mathrm{Mat}_d(R)$. In particular, $B$ is an idempotent lift of $A$.
	
	Let $\widetilde{W}$ denote the image of $B$ in $V$, which is a direct summand of $V$ and hence a free $R$-module. Since its reduction mod $\mathfrak{m}$ is $W$, the rank of $\widetilde{W}$ is $s$, so we have found $\widetilde{W}\in \mathrm{Gr}_s(V)$ mapping to $W$ under the reduction map. In particular $\widetilde{W}/\pi\in \mathrm{Gr}_s(V/\pi)$ maps to $W$, so that the map $\mathrm{Gr}_s(V/\pi)\to \mathrm{Gr}_s(V_k)$ is surjective.
	
	The same idempotent lifting shows that $\mathrm{Gr}_s(V)\to \mathrm{Gr}_s(V/\pi)$ is surjective.
\end{proof}

Given $W\in \mathrm{Gr}_s(V)$, we let $W^\perp=\{v\in V: \omega(v, w)=0\ \forall w\in W\}$, likewise for $\mathrm{Gr}_s(V_k)$.

\begin{lem}
	\label{perpmodm}
	\leavevmode
	\begin{enumerate}[(i)]
		\item If $W\in \mathrm{Gr}_s(V)$, then $W^\perp\in \mathrm{Gr}_{2m-s}(V)$.
		\item If $W\in \mathrm{Gr}_s(V)$, then $(W^\perp)\otimes_Rk\cong (W_k)^\perp$.
	\end{enumerate}
\end{lem}
\begin{proof}
	As in \cite[Lemma 7.2]{AWannals}, it follows from the non-degeneracy of $\omega$ that $W^\perp$ is $\pi$-adically saturated and hence a direct summand of $V$. The rank of $W^\perp$ is then determined by rank-nullity, proving (i).
	
	For (ii), we clearly have $(W^\perp)\otimes_R k\subseteq (W_k)^\perp$. By rank-nullity, both are $k$-vector spaces of dimension $2m-s$, so the result follows.
\end{proof}
Note that $\mathrm{Gr}_s(V_k)$ consists of the $k$-rational points of a $k$-scheme $\mathrm{Gr}_{s, 2m}$.
\begin{lem}
	$\mathrm{Gr}_s(V_k)=\mathrm{Gr}_{s, 2m}(k)$ is an irreducible topological space when viewed as a subspace of $\mathrm{Gr}_{s, 2m}$.
\end{lem}
\begin{proof}
	As $k$ is algebraically closed, $\mathbb{A}^{4m^2}_k(k)$ is dense in the irreducible scheme $\mathbb{A}^{4m^2}_k$, and thus an irreducible topological space. 
	
	As $\mathrm{GL}_{2m}(k)$ is open in $\mathbb{A}^{4m^2}_k(k)$, it is irreducible, and the group $\mathrm{GL}_{2m}(k)$ acts transitively on $\mathrm{Gr}_s(V_k)$. This proves the lemma.
\end{proof}

\begin{lem}
	The map $\perp: \mathrm{Gr}_s(V_k)\to \mathrm{Gr}_{2m-s}(V_k)$ sending $W$ to $W^\perp$ is an isomorphism for each $0\leq s\leq 2m$.
	
	In particular, if $U\subseteq \mathrm{Gr}_m(V_k)$ is a non-empty Zariski-open subset, then there exists some $W\in U\cap U^\perp$.
\end{lem}
\begin{proof}
	The map $W\mapsto \{f\in V_k^*:\ f(w)=0 \ \forall w\in W\}$ provides the usual isomorphism $\mathrm{Gr}_s(V_k)\to \mathrm{Gr}_{2m-s}(V_k^*)$. Since $\omega$ is non-degenerate on $V_k$, $v\mapsto \omega(v, -)$ defines an isomorphism $V_k\cong V_k^*$, inducing by functoriality an isomorphism $\mathrm{Gr}_{2m-s}(V_k^*)\cong \mathrm{Gr}_{2m-s}(V_k)$. The composition of these maps is precisely $W\mapsto W^\perp$.
	
	In particular, if $U\subseteq \mathrm{Gr}_m(V_k)$ is a non-empty Zariski-open subset, then so is $U^\perp$, and $U\cap U^\perp\neq \emptyset$ by the irreducibility of $\mathrm{Gr}_m(V_k)$.
\end{proof}

\begin{prop}
	\label{Weylsimple}
	$\h{A_{m, 2n}}\cong K\langle V\rangle_n$ is simple, i.e. it contains no non-trivial two-sided ideals.
\end{prop}
\begin{proof}
	This is the same argument as in \cite[Theorem 7.3]{AWannals}: if $I\subseteq \h{A_{m, 2n}}$ is a two-sided ideal, then set $\I:=I\cap \h{R_\omega[V]_n}$ and $\I_k=\I\otimes_Rk$. Then as in \cite[Lemma 7.3]{AWannals}, $\I_k$ is an ideal of $k[X_1, \hdots, X_m, Y_1, \hdots, Y_m]$ stable under the divided power derivations $\frac{d_v^k}{k!}$, where $d_v: x\mapsto \frac{vx-xv}{{\pi'}^n}$ for any $v\in V$. Hence $\I_k=0$ or $1\in \I_k$. 
	
	Since $\I$ is $\pi$-adically saturated in $\h{R_\omega[V]_n}$ by construction, it is finitely generated by Corollary \ref{ptfisfp}. Hence it follows from Nakayama's Lemma that if $\I_k=0$, then $\I=0$ and hence $I=0$.
	
	On the other hand, if $1\in \I_k$, then $1\in \I$, since any element in $1+\mathfrak{m}\h{R_\omega[V]_n}$ is a unit. It follows after inverting $\pi$ that either $I=0$ or $1\in I$, so $K\langle V\rangle_n$ is simple.
\end{proof}

\begin{lem}
	\label{Wintdom}
	For any $W\in \mathrm{Gr}_s(V)$, $K\langle W\rangle_n$ is an integral domain.
\end{lem}
\begin{proof}
	It suffices to show that $\h{R_\omega[W]_n}$ is an integral domain.
	
	If $x, y\in \h{R_\omega[W]_n}$ such that $x,y\neq 0$ and $xy=0$, we can assume without loss of generality that $x, y\notin \mathfrak{m}\h{R_\omega[W]_n}$, since $\h{R_\omega[W]_n}$ is $\pi$-torsionfree. 
	
	But now $\h{R_\omega[W]_n}\otimes_Rk\cong \mathrm{Sym}_k W_k$ is an integral domain, providing the desired contradiction.
\end{proof}

The key algebraic input in our proof of Bernstein's inequality (and the only part where we need to add extra work to the discussion in \cite{AWannals}) is the following variant of Noether normalization:

\begin{prop}
	\label{NoetherRpi}
	Let $\A=R/\pi[x_1, \hdots, x_d]$.
	Let $I\subset \A$ be a proper ideal containing $\mathfrak{m}\A$ and let $s$ be the Krull dimension of the $k$-algebra $\A/I$. Then there exists a non-empty Zariski open subset $U\subseteq \mathrm{Gr}_s(k^d)$ and for each $W\in \mathrm{Gr}_s(R^d)$ with $W_k\in U$ a finitely generated ideal $J_W\subseteq I$ such that the following holds:
	\begin{enumerate}[(i)]
		\item $J_W+\mathfrak{m}\A=I$.
		\item $\A/J_W$ is a finitely generated $\mathrm{Sym}_{R/\pi}(W/\pi)$-module.
		\item The natural map
		\begin{equation*}
			\mathrm{Sym}_k(W_k)\to \A/I
		\end{equation*}
		is injective and finite.
	\end{enumerate}
\end{prop}
\begin{proof}
	This is closely modeled on the proof of Noether normalization given e.g. in \cite[section 3.4]{Singular}.
	
	We argue by induction on $d$. If $d=0$, $I$ is a proper ideal of $R/\pi$ containing $\mathfrak{m}$ and hence $I=\mathfrak{m}$. Taking $J=0$, the claim reduces to $R/\pi$ being finitely generated over $R/\pi$, and the map $k\to k$ being injective and finite.
	
	Now let $d\geq 1$. If $I=\mathfrak{m}\A$, then $s=d$, and the claim is again trivial with $J=0$. If $\mathfrak{m}\A\subsetneq I$, we can pick an element $f\in I$ of degree $r\geq 1$ such each non-zero coefficient of $f$ is a unit in $R/\pi$.
	
	Write $f=\sum f_n$ where $f_n$ is homogeneous of total degree $n$. By assumption, the reduction of $f_r$ modulo $\mathfrak{m}$ is a non-zero polynomial over $k$, homogeneous of degree $r$. Now pick any $d\times d$-matrix $M\in \mathrm{GL}_d(R)$ such that $f_r(M_{1d}, \hdots, M_{dd})\notin \mathfrak{m}$ -- this exists as we assume $k$ to be algebraically closed, in particular infinite.
	
	Setting 
	\begin{equation*}
		\begin{pmatrix}
			v_1\\
			v_2\\
			\vdots\\
			v_d
		\end{pmatrix}=M^{-1}\cdot \begin{pmatrix}
		x_1\\
		x_2\\
		\vdots\\
		x_d
		\end{pmatrix},
	\end{equation*}
	we write $V'_M=\sum_{i=1}^{d-1} R\cdot v_i\subseteq \oplus_{i=1}^d R\cdot x_i$. Since the $v_i$ form a basis of $R^d=\oplus_i R\cdot x_i$, $V'_M$ is a free direct summand of $R^d$. In particular, if $W\in \mathrm{Gr}_s(V'_M)$, then $W\in \mathrm{Gr}_s(R^d)$ in a natural way.
	
	It follows that
	\begin{equation*}
		f_r(x_1, \hdots, x_d)=f_r(M_{1d}, \hdots, M_{dd})\cdot v_d^r+(\text{lower terms in $v_d$}),
	\end{equation*}
	so that $f_r=u\cdot v_d^r+(\text{lower terms in $v_d$})$ with $u$ a unit in $R/\pi$, and hence $v_d$ is integral over $R/\pi[v_1, \hdots, v_{d-1}]$ in $\A/(f)$ by construction. Hence the natural map 
	\begin{equation*}
		\phi: R/\pi[v_1, \hdots, v_{d-1}]\to \A/(f)
	\end{equation*}
	is finite. Moreover, the induced map $k[v_1, \hdots, v_{d-1}]\to \A/(f)\otimes_R k$ is injective and finite.
	
	Now let $I'$ be the kernel of the natural map $R/\pi[v_1, \hdots, v_{d-1}]\to \A/(f)\to \A/I$. This is a proper ideal containing $\mathfrak{m}$, so by induction hypothesis, there exists $U_M'\subseteq \mathrm{Gr}_{s'}(k^{d-1})$ and for each $W'$ with $W'_k\in U'$, a finitely generated ideal $J'_W\subseteq R/\pi[v_1, \hdots, v_{d-1}]$ such that $J'_W+\mathfrak{m}=I'$, $R/\pi[v_1, \hdots, v_{d-1}]/J'_W$ is finite over $\mathrm{Sym}_{R/\pi}(W')$ and $\mathrm{Sym}_k(W'_k)\to R/\pi[v_1, \hdots, v_{d-1}]/I'$ is injective and finite, with $s'$ being the Krull dimension of the quotient $R/\pi[v_1, \hdots, v_{d-1}]/J'_W$, which is equal to the Krull dimension of $R/\pi[v_1, \hdots, v_{d-1}]/I'$, as any prime ideal will need to contain the locally nilpotent ideal $\mathfrak{m}/\pi$.
	
	It follows that $\mathrm{Sym}_k(W'_k)\to R/\pi[v_1, \hdots, v_{d-1}]/I'\to \A/I$ is injective and finite, and hence integral. In particular, $s'=s$ by the going-up theorem. Moreover, $\mathrm{Sym}_{R/\pi} (W'/\pi)\to R/\pi[v_1, \hdots, v_{d-1}]/J'_W\to \A/(\phi(J'_W), f)$ is finite -- the first map is finite by assumption on $J'_W$, and the second map is finite since $\phi$ is. 
	
	We now define $U\subseteq \mathrm{Gr}_s(k^d)$ to be the following non-empty Zariski-open subset: $W\in U$ if and only if there exists $M\in \mathrm{GL}_d(R)$ with $f_r(M_{1d}, \hdots, M_{dd})\notin \mathfrak{m}$ such that $W\subseteq V'_M$ and $W\in U'_M$.
	
	Given $W\in \mathrm{Gr}_s(R^d)$ with $W_k\in U$, pick $J'_W$ as above. Lift a finite set of generators of $I/\mathfrak{m}\subseteq k[x_1, \hdots, x_d]$ to $\A$ and adjoin them to $(\phi(J'_W), f)$ to obtain a finitely generated ideal $J_W$ such that $J_W+\mathfrak{m}=I$ with the desired properties.
\end{proof}

\begin{prop}
	\label{Noethersupport}
	Let $\A=R/\pi[x_1, \hdots, x_d]$. Let $\M$ be a finitely generated $\A$-module, $I=\sqrt{\mathrm{Ann}_{\A}(\M)}$, and let $s$ be the Krull dimension of the $k$-algebra $\A/I$. Then there exists a non-empty Zariski open subset $U\subseteq \mathrm{Gr}_s(k^d)$ such that for any $W\in \mathrm{Gr}_s(R^d)$ with $W_k\in U$, $\M$ is almost finitely generated over $\mathrm{Sym}_{R/\pi} (W/\pi)$.
\end{prop}
\begin{proof}
	If $\M=0$, there is nothing to show. If $\M\neq 0$, we know that $I$ is a proper ideal of $\A$ containing $\mathfrak{m}$, so we can consider the non-empty Zariski-open subset $U\subseteq \mathrm{Gr}_s(k^d)$ provided by Proposition \ref{NoetherRpi}.
	
	If $W\in \mathrm{Gr}_s(R^d)$ with $W_k\in U$, let $J_W\subseteq I$ be a finitely generated ideal of $\A$ as in Proposition \ref{NoetherRpi}. Since $J_W$ is finitely generated, there exists some $r$ such that $(J_W)^r\cdot \M=0$. We thus have a finite filtration on $\M$ with
	\begin{equation*}
		\M_i=\{m\in \M: (J_W)^i\cdot m=0\}
	\end{equation*}	
	such that $\M_i/\M_{i-1}$ is an almost finitely generated left $\A/J_W$-module.
	
	But by Proposition \ref{NoetherRpi}, $\M_i/\M_{i-1}$ is then an almost finitely generated $\mathrm{Sym}_{R/\pi}(W/\pi)$-module, and thus $\M$ is almost finitely generated over $\mathrm{Sym}_{R/\pi}(W/\pi)$.
\end{proof}

\begin{thm}
	\label{Bernsteinsupport}
	Let $\M$ be a finitely generated $\pi$-torsionfree $\h{R_\omega[V]_n}$-module. Then $\mathrm{Ext}^j_{\h{R_\omega[V]_n}}(\M, \h{R_\omega[V]_n})\otimes_R K= 0$ for all $j>m$. 
\end{thm}

\begin{proof}
	Fix $j>m$. We need to show that $\N:=\mathrm{Ext}^j_{\h{R_\omega[V]_n}}(\M, \h{R_\omega[V]_n})$ is $\pi$-torsion. Let us assume that $\N$ is not almost zero (otherwise, there is nothing to show).
	
	Note by Lemma \ref{Extfg}.(ii) that $\N$ is a finitely generated $\h{R_\omega[V]_n}$-module, so by Nakayama's Lemma, $\N/\N\mathfrak{m}\neq 0$. Since 
	\begin{equation*}
		\N/\N\mathfrak{m}\subseteq \mathrm{Ext}^j_{\h{R_\omega[V]_n}}(\M, \h{R_\omega[V]_n}/\mathfrak{m})\cong \mathrm{Ext}^j_{\h{R_\omega[V]_n}/\mathfrak{m}}(\M/\mathfrak{m}, \h{R_\omega[V]_n}/\mathfrak{m})
	\end{equation*}
	it follows from the regularity of $\h{R_\omega[V]_n}/\mathfrak{m}\cong k[X_1, \hdots, X_m, Y_1, \hdots, Y_m]$ that the Krull dimension $s$ of $\h{R_\omega[V]_n}/I$ for $I=\sqrt{\mathrm{Ann}(\N/\N\mathfrak{m})}$ is at most $2m-j$, which is strictly less than $m$.
	
	But now $I=\sqrt{\mathrm{Ann}_{\h{R_\omega[V]_n}/\pi}(\N/\N\pi)}$ and $\h{R_\omega[V]_n}/\pi\cong R/\pi[X_1, \hdots, X_m,Y_1, \hdots, Y_m]$, so Proposition \ref{Noethersupport} provides us with a non-empty, Zariski open subset $U\subseteq \mathrm{Gr}_s(V_k)$ such that $\N/\N\pi$ is almost finitely generated over $\mathrm{Sym}_{R/\pi}(W/\pi)$ for all $W\in \mathrm{Gr}_s(V)$ with $W_k\in U$. A fortiori, there exists a non-empty Zariski open subset $U'\subseteq \mathrm{Gr}_m(V_k)$ such that $\N/\N\pi$ is almost finitely generated over $\mathrm{Sym}_{R/\pi}(W/\pi)$ for any $W\in \mathrm{Gr}_m(V)$ with $W_k\in U'$, as $s<m$.
	
	Since $\mathrm{Gr}_m(V_k)$ is irreducible and $\perp$ is a homeomorphism, we can find $\overline{W}\in U'\cap {U'}^\perp$. Since $\mathrm{Gr}_m(V)\to \mathrm{Gr}_m(V_k)$ is surjective, we then have some $W\in \mathrm{Gr}_m(V)$ such that $W_k=\overline{W}$, and hence $\N/\N\pi$ is almost finitely generated over $\mathrm{Sym}(W/\pi)$ and over $\mathrm{Sym}(W^\perp/\pi)$ (this uses that $\perp$ commutes with reduction modulo $\mathfrak{m}$, Lemma \ref{perpmodm}).
	
	As $\N$ is finitely generated over $\h{R_\omega[V]_n}$, it is $\pi$-adically separated by Lemma \ref{almostAR}.(ii). Hence $\N$ is almost finitely generated over $R_\omega\langle W\rangle_n$ and over $R_\omega\langle W^\perp\rangle_n$ by Lemma \ref{afgcompletion}.(i). In particular, $\N\otimes_RK=\mathrm{Ext}^j_{K\langle V\rangle_n}(\M\otimes_RK, K\langle V\rangle_n)$ is finitely generated over $K\langle W\rangle_n$ and $K\langle W^\perp\rangle_n$. 
		
	From here on, the proof proceeds as in \cite{AWannals}: let $m_1, \hdots, m_r\in \N\otimes_RK$ be a finite set of generators of $\N\otimes_RK$ over $K\langle W^\perp\rangle_n$. Since $W$ and $W^\perp$ commute with each other, we have
	\begin{equation*}
		\mathrm{Ann}_{K\langle W\rangle_n}=\cap \mathrm{Ann}_{K\langle W\rangle_n} (m_i).
	\end{equation*}
	Each $\mathrm{Ann}(m_i)$ is non-zero, else we could find an $R_\omega\langle W\rangle_n$-submodule of $\N$ which is free of rank $1$ -- but $\N$ is almost finitely generated over $R_\omega\langle Z\rangle_n$ for some $Z\in \mathrm{Gr}_s(W)$, $s<m$, and $R_\omega\langle W\rangle_n$ is not almost finitely generated over $R_\omega\langle Z\rangle_n$. So it follows from \cite[Theorem 2.1.15]{Robson} that $\mathrm{Ann}_{K\langle W\rangle_n}(\N_K)\neq 0$, since $K\langle W\rangle_n$ is an integral domain by Lemma \ref{Wintdom}. But then $\mathrm{Ann}_{K\langle V\rangle_n}(\N_K)\neq 0$ is a non-trivial two-sided ideal, so by simplicity (Proposition \ref{Weylsimple}), $\mathrm{Ann}_{K\langle V\rangle_n}(\N_K)=K\langle V\rangle_n$ and $\N_K=0$, as required.
\end{proof}

\begin{proof}[Proof of Theorem \ref{Bernstein}]
	Let $A=\h{A_{m, 2n}}$ and let $K'$ be a complete field extension of $K$ with algebraically closed residue field. We set $A'=A\h{\otimes}_K K'$, the corresponding Tate-Weyl algebra over $K'$. Let $R'$ denote the valuation ring of $K'$.
	
	We first show that $\mathrm{gl.dim.}(A')\leq m$.
	
	Let $M$ be a finitely generated $A'$-module, and let $\M\subseteq M$ be a finitely generated $\pi$-torsionfree $\h{R'_\omega[V]_n}$-module such that $\M\otimes_{R'}K'\cong M$. Then by Theorem \ref{Bernsteinsupport}, $\mathrm{Ext}^j_{A'}(M, A')\cong \mathrm{Ext}^j_{\h{R'_\omega[V]_n}}(\M, \h{R'_\omega[V]_n})\otimes_{R'}K'=0$ for all $j>m$, so the global dimension of $A'$ is at most $m$.
	
	Now Lemma \ref{injbasechange} implies that $\mathrm{gl.dim.}(A)\leq m$.
	
	To show that the global dimension is not smaller than $m$, we can proceed analogously to subsection 5.1:
	
	Taking $T=K\langle x_1, \hdots, x_m\rangle$, $L=\mathrm{Der}_K(T)$, the Spencer resolution from \cite[subsection 6.3]{SixOp} produces a free resolution
	\begin{equation*}
		\hdots\to A\otimes_T\wedge^i L\to \hdots \to A\otimes_T\wedge^0 L\to T\to 0
	\end{equation*}
	of $T$ as a left $A$-module, giving us
	\begin{equation*}
		\mathrm{Ext}^i_{A}(T, A)=\begin{cases}
			\Omega^m\ \text{if } i=m\\
			0\ \text{otherwise},
		\end{cases}
	\end{equation*}
	where $\Omega_m=\wedge^m \Omega^1_X(X)$ for $X=\mathrm{Sp}T$ the $m$-dimensional polydisc.
	
	Thus $\h{A_{m, 2n}}$ has global dimension exactly $m$. As we can replace $\pi$ by a square root of $\pi$ after possibly a finite field extension by Lemma \ref{injbasechange}, this yields $\mathrm{gl.dim.}\h{A_{m, n}}=m$ for all $n\geq 1$.
	
	In particular, if $M$ is a finitely generated left $\h{A_{m, n}}$-module with
	\begin{equation*}
		\mathrm{Ext}^j_{\h{A_{m, n}}}(M, \h{A_{m,n}})=0
	\end{equation*}
	for all $j\leq m$, then $\mathrm{R}\mathrm{Hom}_{\h{A_{m,n}}}(M, \h{A_{m,n}})=0$, so 
	\begin{equation*}
		M\cong \mathrm{R}\mathrm{Hom}_{\h{A_{m,n}}}(\mathrm{R}\mathrm{Hom}_{\h{A_{m,n}}}(M, \h{A_{m,n}}), \h{A_{m,n}})=0.
	\end{equation*}
	This finishes the proof.
\end{proof}

\end{document}